\documentclass[11pt,reqno]{amsart}

\title[]{Electroconvection in a Magnetic Field}

\author[E. Abdo]{Elie Abdo}
\address[E. Abdo]
{	Department of Mathematics \\
     American University of Beirut \\
	Beirut 1107-2020\\Lebanon.} \email{ea94@aub.edu.lb}

\author[P. Constantin]{Peter Constantin}
\address[P. Constantin]
{	Department of Mathematics \\
     Princeton University  \\
	Princeton, NJ 08544, USA.} \email{const@math.princeton.edu}

\author[M. Ignatova]{Mihaela Ignatova}
\address[M. Ignatova]
{	Department of Mathematics \\
     Temple University  \\
	Philadelphia, PA 19122, USA.} \email{ignatova@temple.edu}

\author[Q. Lin]{Quyuan Lin}
\address[Q. Lin]
{	School of Mathematical and Statistical Sciences \\
Clemson University\\
Clemson, SC 29634, USA.} \email{quyuanl@clemson.edu}

\usepackage[margin=1in]{geometry}
\usepackage{amsmath, amsthm, amssymb}
\usepackage{times}
\usepackage{comment}
\numberwithin{equation}{section}
\usepackage{color}
\usepackage[colorlinks=true, pdfstartview=FitV, linkcolor=blue, citecolor=blue, urlcolor=blue]{hyperref}
\usepackage{hyperref}
\usepackage{enumerate}
\newcommand{\pa}{\partial}

\newcommand{\la}{\label}
\newcommand{\fr}{\frac}
\newcommand{\na}{\nabla}
\newcommand{\be}{\begin{equation}}
\newcommand{\ee}{\end{equation}}
\newcommand{\ba}{\begin{array}{l}}
\newcommand{\ea}{\end{array}}
\newcommand{\eps}{\epsilon}

\newcommand{\beg}{\begin}

\renewcommand{\div}{{\mbox{div}\,}}
\newcommand{\D}{\Delta}
\renewcommand{\l}{\Lambda}

\usepackage{xcolor}
\usepackage{graphicx}
\usepackage[]{amsfonts}
\usepackage{amsthm}
\usepackage{amsmath}
\usepackage {amssymb}
\usepackage{bm}
\usepackage{mathrsfs}
\usepackage{setspace}
\usepackage{amsmath}
\usepackage{cite}

\usepackage{MnSymbol,wasysym}

\newcommand{\N}{\mathbb N}

\def\RR{{\mathbb R}}
\def\TT{{\mathbb T}}
\def\NN{\mathbb N}
\def\PP{\mathbb P}
\renewcommand{\div}{\na\cdot\,}

\date{\today}

\begin{document}

\begin{abstract}
 Electroconvection in a porous medium under a strong transversal magnetic field is described by an active scalar equation for the charge density. The equation has global weak solutions with
 $L^{\infty}$ data. We show that for strong enough magnetic fields, $L^{\infty}$-small solutions are smooth globally in time and they obey surface quasigeostrophic equations in the limit of infinite magnetic field strength.
 \end{abstract}

\noindent

\maketitle

\section{Introduction}
 Electroconvection, the flow of charges in fluids, is characterized by the fact that the charge density and  the fluid's velocity and pressure are directly coupled. The charges in the solvent exert forces on the fluid. The fluid responds to these forces while convecting the charges.  This results in a nonlinearly coupled system of 
 equations for the charge density and the fluid's velocity and pressure.  The subject belongs to a large class of electrohydrodynamic problems of broad scientific interest. 
Experimental and theoretical studies \cite{morris1991patterns,daya1997electroconvection,daya1999electrically,daya1998annular} on  smectic films revealed complex dynamical behavior in rotating, two-dimensional conductive fluids subjected to three-dimensional electrostatic forces. In experimental studies, applied transversal magnetic fields  significantly impact the dynamics of the electroconvection \cite{gleeson1996onset, ahmed2023effect, huh1999formation, huh2007electrohydrodynamic, huh2015multiplicative}. In particular, when a magnetic field is applied,  the critical voltage required for instability increases, leading to enhanced system stability \cite{ahmed2023effect}. 
Magnetic fields applied to ordered electrically sensitive fluid, such as liquid crystals in porous media, have been documented to produce modifications in phase transitions, changes in orientational order, alterations in the elastic properties and in the behavior of the director field \cite{dadmun1993nematic, sinha1998dielectric, goldburg1995behavior}. 

 In this paper we are concerned with mathematical properties of electroconvection in a magnetic field.
  We show that in the limit of strong transversal magnetic field, the solutions of equations of electroconvection  in porous media converge to solutions of the surface quasigeostrophic equation (SQG).

The charge density $q(x,t)$ of electroconvection in porous media is an active scalar. Active scalars \cite{constantin1994geometric, pierrehumbert1994spectra} are transported by incompressible velocities they create by means of a time-independent equation. This widely studied class of equations  includes the two-dimensional incompressible vorticity equation, the surface quasigeostrophic equation SQG and generalized g-SQG equations interpolating between them.
The SQG  active scalar arose in geophysics \cite{held1995surface} as a model of the large scale mid-latitude surface temperature evolution in quasigeostrophic  flow. Electroconvection of charges in porous media under strong magnetic fields and atmosphere-ocean thermal dynamics are such disparate physical systems, and yet, mathematically, they turn out to be related.
The inviscid SQG equation is studied in the context of singularity formation in fluids \cite{constantin1994geometric,cordoba1998nonexistence,constantin1994formation}.  Global regularity from arbitrary smooth data is not known. The equation is ill posed in spaces of low regularity \cite{cordoba2024instantaneous,cordoba2024non}. Global weak solutions were obtained in \cite{resnick1996dynamical} (see also \cite{constantin2018local} for bounded domains). Well-posedness and loss of regularity results have been obtained for g-SQG in the case the initial data is a patch (that is a step function) with smooth boundary in a half-space  \cite{kiselev2017local, kiselev2016finite, zlatos2023local, gancedo2021local} and fronts \cite{berti2025paralinearization, gancedo2022well}. Solutions that exist for all time have been constructed. They include steady radial solutions, time-periodic rotating solutions \cite{castro2020global} and quasiperiodic solutions (see the monograph \cite{serrano2023quasiperiodic} and references therein). The critical dissipative SQG equation has global smooth solutions from arbitrary data \cite{kiselev2007global,caffarelli2010drift}. Critical dissipative SQG is $L^{\infty}$ critical. Small $L^{\infty}$ data lead to stability and global existence \cite{constantin2001critical}; for large data, an initial data dependent modulus of continuity \cite{kiselev2007global} or a $C^{\alpha}$ norm \cite{constantin2015long} are non-increasing in time. These quantities depend on the size of the initial data in $L^{\infty}$ and are highly sensitive to the structure of the equation: finding them relies on the linear nonlocal equation of state of the velocity in terms of the scalar. This is a source of difficulty in analyzing nonlinear nonlocal top-order perturbations of the critical SQG, such as we encounter in the present work. Once these quantities are finite, higher regularity follows \cite{constantin2008regularity}.

We describe now the electroconvection setting we are discussing in this paper.  We consider an incompressible fluid occupying a very thin region in space,  represented as the two dimensional plane $\Omega= \RR^2 =\{(x,y,z)\left| \right.\, z=0\}\subset \RR^3$. We denote by $q$ a charge density in  $\Omega$, denote by $\tilde{E}=(\tilde E_1,\tilde E_2,\tilde E_3)$ and $\tilde{B}=(\tilde B_1,\tilde B_2,\tilde B_3)$ three-dimensional electric and magnetic fields, and by $\rho$ the total charge density $\rho = 2q\delta_{\Omega}$ as a distribution in the whole space. We also denote by $E=(E_1, E_2)$ the restriction of $(\tilde E_1,\tilde E_2)$ to $\Omega$. The three-dimensional electric field $\tilde{E}$  obeys Gauss' law 
\begin{equation} 
\na_3 \cdot \tilde{E} = \frac1{\epsilon_0}\rho = \frac1{\epsilon_0}2q\delta_{\Omega},
\end{equation}  where $\epsilon_0$ is the vacuum permittivity and we denote by $\na_3$ the three dimensional gradient. The electric field is given by an electrostatic potential $\Phi$,
\begin{equation} 
\tilde{E} = - \na_3 \Phi. 
\end{equation}  Gauss' law gives the Poisson equation obeyed by $\Phi$, 
\begin{equation} 
-\Delta_3 \Phi = \frac1{\epsilon_0}2q \delta_{\Omega}.
\end{equation}   
We discuss the case when $q$ and $\Phi$ are $2\pi$ periodic functions of $x$ and $y$ and $\Phi$ vanishes as $|z|\to\infty$.  The potential $\Phi $ is  then given by
\begin{equation} 
\Phi (x,y,z) = \begin{cases}
\frac1{\epsilon_0}e^{-z \Lambda} \Lambda^{-1} q, \hspace{1cm} z > 0,
\\ \frac1{\epsilon_0}e^{z \Lambda} \Lambda^{-1} q, \hspace{1cm} z < 0,
\end{cases}
\end{equation}  where $\Lambda = \sqrt{-\Delta}$ is the square root of the two dimensional  Laplacian with $2\pi$ periodic boundary conditions. Therefore, $E:=(\tilde{E}_1, \tilde{E}_2)|_{z=0}$ can be expressed as
\begin{equation} 
 E = -\frac1{\epsilon_0}\na \l^{-1} q 
 \la{eleq}
\end{equation} 
where $\na$ is the two dimensional gradient. We denote by $R$ the Riesz transforms 
\be
R = \na\l^{-1}
\la{riesz}
\ee
and thus the electric field  restricted to the domain occupied by fluid is
\be
E = - \fr{1}{\eps_0}Rq.
\la{elecRq}
\ee
The magnetic field generated by the electric field is negligible. We consider an imposed external magnetic field $\tilde{B}=(0,0,B)$ where $B\geq 0$ is a constant independent of both time and space.  Notice that this external magnetic field $\tilde{B}$ satisfies the divergence-free condition, has zero curl, and is time-independent.

The total current density is given by the sum of the advective current $(u_1, u_2, 0)q$ and the Ohmic conduction current $\sigma \tilde{E}'$,  where $\sigma>0$ is  the conductivity of the medium and $\tilde{E}'$ is the electric field experienced by the fluid element in its rest frame \cite{levan2025foundations}. When the fluid is moving with
respect to the external magnetic field at the velocity $u$, using the Lorentz transformation we have
\begin{equation*}
    \tilde{E}'= \tilde{E} + (u_1,u_2,0)\times \tilde{B}.
\end{equation*} 
Therefore, 
\begin{equation}
    \tilde{j} = (u_1, u_2, 0)q + \sigma (\tilde E + (u_1,u_2,0)\times \tilde{B}).
\end{equation}

The restriction of $\tilde{j}$ to the surface is denoted by $j=(\tilde{j}_1,\tilde{j}_2)|_{z=0}$ and is therefore
\begin{equation}
    j = uq + \sigma (E - Bu^\perp ),
    \la{curruEB}
\end{equation}
where $u^\perp= (-u_2, u_1)$. The surface charge density $q$ obeys the conservation 
\be
\pa_t q + \div j =0
\la{consq}
\ee
which is, in view of \eqref{curruEB},
\begin{equation}
   \partial_t q + u\cdot \nabla q + \sigma (\frac1{\epsilon_0}\Lambda q - B \nabla\cdot u^\perp)=  0.
\la{qtransp}
\end{equation} 
The fluid's divergence-free velocity $u$ obeys an equation forced by the electrostatic force $qE$ plus the Lorentz force
$(\tilde j\times \tilde B)$ in the plane, which is computed as
\be
(\tilde j\times \tilde B)_{\left |\right. \RR^2} = \fr{\sigma}{\eps_0} BR^{\perp} q - Bq u^{\perp} -\sigma B^2 u.
\la{lorentz}
\ee
To ease the notation we take $\sigma= \eps_0 =1$. The total force exerted on the fluid is given by
\be
F =  BR^{\perp} q - qRq  - Bqu^{\perp} - B^2 u.
\la{force}
\ee
There are two important velocity contributions in $F$ due to the imposed external magnetic field. 
The first one, $-Bqu^{\perp}$, yields a rotation effect similar to a Coriolis force with frequency $Bq$.   The 
second one $- B^2 u$, is  a strong damping or friction effect.
Electroconvection  equations couple \eqref{qtransp} to a momentum equation driven by $F$. For instance, the Navier-Stokes equations forced by $F$ are
\be
\pa_t u + u\cdot\na u -\nu\D u +  B^2 u + Bqu^{\perp} + \na p = BR^{\perp}q -qRq, \quad \na\cdot u =0,
\la{nse}
\ee
where $\nu>0$ is the kinematic viscosity. In the absence of an applied transversal magnetic field, the Navier-Stokes electroconvection  system, i.e the system \eqref{qtransp}, \eqref{nse} with $B=0$, was shown to have global solutions in \cite{constantin2017some}. The long time dynamics were described in \cite{abdo2021long}  in $\TT^2$ and in \cite{abdo2025long} in bounded domains, both with,  and without  time independent body forces in the fluid.

 If we replace viscous friction by Darcy's law and neglect inertial time dependence, we arrive at the equations
\be
\mu\, u + \na p = F, \quad \na\cdot u =0,
\la{darcy}
\ee
where $\mu>0$ is porosity. This is the fluid equation we consider in this paper. 
Taking $\mu=1$ and retaining $B$ as the only  (large) parameter, \eqref{darcy} becomes the time independent law
\be
(1+B^2)u + Bq u^{\perp} + \na p = BR^{\perp}q - qRq, \quad \na\cdot u = 0.
\la{uqp}
\ee 
Together with the charge density equation \eqref{qtransp}, this gives rise to the following system,
\beg{align}  
& \partial_t q + u\cdot \nabla q + \Lambda q - B \nabla\cdot u^\perp=  0, \label{eqn:q}
\\
&  (1+ B^2) u + Bqu^{\perp} + \nabla p = BR^{\perp} q-qRq, \label{eqn:u}
\\
&\na \cdot u =0 \label{eqn:div}.
\end{align}
The unknowns are $q$, $u$, $p$. Using 
$-B\nabla\cdot u^{\perp} =B\omega$, taking the two dimensional curl of \eqref{eqn:u},  and replacing in \eqref{eqn:q}, the latter becomes
\be
\pa_t q + \fr{B}{1+B^2}(R^{\perp}q\cdot\na q) + \fr{1}{1+B^2}u\cdot\na q +\fr{1}{1+B^2}\l q =0.
\la{decm}
\ee
In the absence of an external magnetic field, i.e. when $B=0$, the system \eqref{eqn:q}-\eqref{eqn:div}   describes  electroconvection in a porous medium. The global existence for large data of smooth solutions of this system is a challenging open problem.  Global regularity for small initial data in Besov spaces smaller than $L^{\infty}$ was obtained in \cite{abdo2022electroconvection}.

We consider three time scales, $t=t_{lab}$ the laboratory time scale,
$t_1$, the magnetic gyration time scale or first magnetic time scale, and $t_2 $, the magnetic friction time scale or second magnetic time scale. They are related by
\be
t_1 = \fr{B}{1+B^2} t \quad \text{and}\quad t_2 = \fr{1}{B}t_1 = \fr{1}{1+B^2} t.
\la{times}
\ee
In laboratory time scale, the equation for the unknown $q(x,t)$ is \eqref{decm}, that is
\be
\pa_t q + \fr{B}{1+B^2}(R^{\perp}q\cdot\na q) + \fr{1}{1+B^2}\left (u[q]\cdot\na +\l\right)q =0,
\la{decmlab}
\ee
with $u[q] = u$  the unique solution of
\be
(1+B^2)\left (u + \fr{B}{1+B^2}\mathbb P (qu^{\perp})\right ) =  BR^{\perp} q - \mathbb P(qRq)
\la{u[q]eq}
\ee
where $\PP = \mathbb I + \na (-\Delta)^{-1} \na \cdot = \mathbb I  + R\otimes R$ is the Leray projector on divergence free vector fields. 
The time independent equation of state $u[q]$  and its properties are described in detail below in Section \ref{ES}. 

In the first (gyration) magnetic time scale, the equation for $q_1(x, t_1) = q(x, \fr{1+B^2}{B}t_1)$ where $q$ solves \eqref{decmlab} becomes
\be
\pa_{t_1} q_1 + R^{\perp}q_1\cdot\na q_1 + \fr{1}{B}\left (u[q_1]\cdot\na +\l\right)q_1 = 0.
\la{decmtm}
\ee
The initial data are the same, $q(x, 0) = q_1(x, 0)$, and there is no rescaling of space variable $x$.

 In the second magnetic time scale, the equation for the boosted field
\be
Q(x, t_2) = B q_1(x, Bt_2) = Bq(x, (1+B^2) t_2)
\la{Qq}
\ee 
becomes
\be
\pa_{t_2}Q +  R^{\perp}Q\cdot\na Q + \left (u\left[\fr{Q}{B}\right]\cdot\na +\l\right)Q = 0.
\la{rdecmtm}
\ee
The initial data is $Q(x,0) = Bq(x,0)$ and there is no rescaling of space variable.  Note that $Q= Bq$ is the frequency of rotation produced by the force $F$.

For fixed $B$, the three 
equations \eqref{decmlab}, \eqref{decmtm}, \eqref{rdecmtm} are versions of the same equation, and they are entirely equivalent. We refer to these equations as Darcy law electroconvection in a magnetic field, DECM equations.

The main results of this work are as follows. We prove first that  weak solutions of \eqref{decm} with initial data in $L^{\infty}$ exist globally (Theorem \ref{wwweak}). 
As $B\to\infty$, in laboratory time scale, the equation converges to $\pa_t q=0$. It is only in the first magnetic time scale that nontrivial dynamics arise. The solutions $q_1$ of  \eqref{decmtm} converge to solutions of inviscid SQG. We prove (Theorem \ref{invsqglim})
\beg{thm}\la{wsqglim}
Any family of weak solutions $q^B$of \eqref{decm}
 has a subsequence such that $q^B(x, \fr{1+B^2}{B}t)$ converges weakly as $B\to \infty $ to a weak solution $q(x,t)$ of the inviscid SQG 
 \[
 \pa_t q + (R^{\perp}q)\cdot\na q = 0,
 \]
  in $L^2(0, T; L^2)$.
\end{thm}

The dissipative nature of the DECM equation emerges in longer time scales, and the boosted field $Q$ converges in the second magnetic time scale to solutions of critical dissipative SQG. We prove global regularity of the boosted field $Q$ (Theorem \ref{globQ}) for small $L^{\infty}$ initial data and large enough $B$.
In laboratory time scale, the regularity result is (Remark \ref{globlab})
\beg{thm}\la{globregq}
There exist $c>0$ and $C>0$ such that, if $q_0$ obeys $B\|q_0\|_{L^{\infty}} \le c$ and $\|q_0\|_{H^3}\le C(1+B^2)$, then the solution of \eqref{decm} with initial data $q_0$ exists globally, is unique, and satisfies
\[
\|q(\cdot,t)\|_{H^3} \le \|q_0\|_{H^3}e^{-\fr{t}{4(1+B^2)}}
\]
for all $t\ge 0$.
\end{thm}
The regularity in $H^3$ implies by parabolic estimates $C^{\infty}$ regularity, a fact that  follows by well known methods and is not pursued in this paper.  

We prove convergence to dissipative critical SQG without assuming the condition of small 
initial $L^{\infty}$ norm of $q$, but we do need to assume enough regularity to allow for absolute continuity of the $L^2$ norm of the DECM solution. Of course, this regularity is guaranteed for short time or for arbitrary time if the initial $L^{\infty}$ norm is small. The result we prove is (Theorem \ref{critsqgthm} and  Remark \ref{critlab})
\beg{thm}
Let $q$ be a strong solution of \eqref{decm} on $[0, (1+B^2)T]$ with initial data $q_0=\frac{1}{B} Q_0$, with fixed $Q_0\in H^3$.  There exists a constant $C$ depending only on $\|Q_0\|_{H^3}$ and $T$ such that 
 \[
\sup_{t\in[0,T]}\|Bq(\cdot, (1+ B^2)t)-\overline{Q}(t)\|_{L^2} \le\fr{C}{1 +B^2}
\]
where $\overline Q$ is the global smooth solution of critical SQG 
\[
\pa_t Q + (R^{\perp}Q)\cdot\na Q + \l Q = 0
\]
with initial data $Q_0$.
\end{thm}
The constant $C$ above is in fact finite as soon as $Q_0\in C^{\alpha}$, for any $\alpha>0$. The result 
holds for fixed $Q_0$ while $B$ varies, and not for fixed $q_0$. If we fix $q_0$ and consider the family of initial data $Bq_0$ for the critical SQG, then the constant $C$ depends badly on $B$.

The paper is organized as follows. In Section \ref{ES} we discuss in detail the equation of state for the velocity.  We show that solving for the velocity in terms of $q$ yields a unique solution
\[
u[q] = R^{\perp}\left( w[q] + \fr{B q}{1+B^2}\right),
\]
where $w$ is small in $H^3$ compared to $\|q\|_{H^3}^2$ for large  $B$ (Proposition \ref{goodwQbound}).
This is a crucial ingredient in the proof of global regularity. 
In Section \ref{S2} we prove the global existence of weak solutions by introducing a judicious approximation.
Section \ref{S3} is where we prove convergence to inviscid SQG equations in the first magnetic time scale.
In Section \ref{S4} we prove the global regularity and dissipative SQG limit in the second magnetic time scale. We provide the proof of Proposition \ref{goodwQbound} in Appendix~\ref{A1} and the proof of global regularity of the approximation in Appendix~\ref{A2}.

\section{Equation of State for the Velocity in Terms of the Charge Density}\la{ES}

Taking the curl of \eqref{u[q]eq} we obtain the time-independent  equation for the vorticity,
\be
\omega = \fr{1}{1+B^2}\left(R^{\perp}q\cdot\na q\right) - \fr{B}{1+B^2}\left (u\cdot\na q + \l q\right),
\la{omega}
\ee
where 
\be
\omega = \na^{\perp}\cdot u = \Delta\psi, \quad \na\cdot u =0.
\la{curl}
\ee
Writing $\omega = -\l v$ where $v=\l \psi$ with $\psi$ the stream function, 
we have from \eqref{omega} after applying $-\l^{-1}$ to both sides 
\be
v -\fr{B}{1+B^2}R\cdot(qR^{\perp}v) = \fr{B}{1+B^2}\left(q - \fr{1}{B} R\cdot(qR^{\perp} q)\right),
\la{veq}
\ee  
and
\be
u = R^{\perp} v.
\la{uv}
\ee
We denote by $L_q$ the operator
\be
L_q (f) = R\cdot(qR^{\perp}f)
\la{Lq}
\ee
Both $q$ and $f$ are scalar.  Because $R\cdot R^{\perp}=0$, this is a commutator,
\be
L_q(f) = [R,q]\cdot R^{\perp} f.
\la{Lqcommut}
\ee
Now \eqref{veq} can be written as
\be
v- \fr{B}{1+B^2}L_q v = \fr{B}{1+B^2}\left (q -\fr{1}{B}L_q q\right).
\la{vLq}
\ee

The operator $\mathbb I - \fr{B}{1+B^2}L_q$ is bounded in $L^2$ (if $q\in L^{\infty}$) and unconditionally invertible. Because $q$ is real, the operator $L_q$ is bounded anti-selfadjoint in $L^2$ so its spectrum lies in a segment of $i\mathbb R$. The inverse of  $\mathbb I - \fr{B}{1+B^2}L_q$ can be defined by familiar functional calculus, and has norm less than 1,  no matter how large is the norm of $\|q\|_{L^{\infty}}$. We recall that under DECM equations this norm is not growing in time.

Introducing
\be
w = v-\fr{B}{1+B^2}q
\la{wq}
\ee
the equation  \eqref{vLq} becomes
\be
w- \fr{B}{1+B^2}L_q w = -\fr{1}{(1+B^2)^2} L_q q,
\la{weqq}
\ee
that is,
\be
w[q] = -\fr{1}{(1+B^2)^2}\left (\mathbb I - \fr{BL_q}{1+B^2}\right)^{-1}(L_q q)
\la{w[q]}
\ee
and \eqref{uv} becomes
\be
u[q] = R^{\perp}\left( w[q] + \fr{B q}{1+B^2}\right).
\la{uw}
\ee
The left hand side of \eqref{u[q]eq} defines an operator  $T_q$  in $H$,  the space of periodic divergence-free vector fields in $L^2$:
\be
T_q(u) = (1+B^2)\left (u + \fr{B}{1+B^2}\mathbb P (qu^{\perp})\right ).
\la{Tq}
\ee
Using \eqref{uv} we express in terms of $v$,
\be
T_q(u) = (1+B^2)\left (R^{\perp}v - \fr{B}{1+B^2}\mathbb P(qRv)\right).
\la{Tquv}
\ee
Now, in view of 
\be
 (\mathbb P f)_i = \left(\delta_{ij} +R_iR_j\right)f_j
 \la{mathbbP}
 \ee
\be
 L_q(v) = - R^{\perp}\cdot(qRv),
 \la{lqswitch}
 \ee
 and 
 \be
 R_iR_i= -\mathbb I,
 \la{Rsquare}
 \ee
by applying from the left  $R^{\perp}_i = \epsilon_{ji}R_j$  to \eqref{Tquv} (where $\epsilon_{ji} =0$ if $j=i$, and it is the signature of the permutation $(1,2)\mapsto (j,i)$ if $j\neq i$),  we obtain 
\be
R^{\perp}\cdot T_q(u) = -(1+B^2)\left(\mathbb I -\fr{B}{1+B^2} L_q\right )(v). 
\la{TqLq}
\ee
This shows that $u$ determined by inverting $T_q$ is the same as $u[q]$.
Indeed, using \eqref{vLq} to solve for $v$, we deduce

\be
R^{\perp}\cdot T_q(u) = -(1+B^2)\fr{B}{1+B^2}\left (q -\fr{1}{B}L_q q\right) = - Bq + L_q q.
\la{lqsolved}
\ee
Now, because $T_q(u)$ is divergence-free, we have
\be
T_q(u) = -R^{\perp}(R^{\perp}\cdot T_q(u)),
\la{TquRperp}
\ee
and using \eqref{lqsolved} we deduce
\be
T_q(u) = BR^{\perp} q - R^{\perp}(L_q q).
\la{Tqsolved}
\ee
Because $\mathbb P(qRq)$ is divergence  free we have that
\be
\mathbb P(qRq) = -R^{\perp}\left (R^{\perp}\cdot \mathbb P(qRq)\right ),
\la{qrqrperp}
\ee
and using \eqref{lqswitch} we see that
\be
\mathbb P(qRq) =  R^{\perp}(L_q q)
\la{pqrqlqq}
\ee
and thus, from \eqref{Tqsolved}
that 
\be
T_q(u) = BR^{\perp} q - \mathbb P(qRq).
\la{tqueq}
\ee
This shows that $u=R^{\perp}v$ solves \eqref{tqueq}. The solution is unique, and this establishes the converse implication as well. We have thus, 
\beg{prop}\la{eqstate} Let $q\in L^{\infty}$. The unique solution in $H$ of \eqref{u[q]eq},
\be
u[q] = T_q^{-1}\left (BR^{\perp}q - \mathbb P(qRq)\right)
\la{u[q]Tq}
\ee
is given by \eqref{uw} with $w[q]$ solving \eqref{weqq}  given by \eqref{w[q]}.
\end{prop}

Note that $L_q(f)$ is bilinear, in particular $BL_q= L_{Bq} $.  In view of \eqref{w[q]} and \eqref{uw}  we have
\be
u[q] = \fr{B}{1+B^2}R^{\perp}\left( q - \fr{1}{B^2}\sum_{n=0}^{\infty}\left(\fr{BL_q}{1+B^2}\right)^{n+1} q\right).
\la{uq}
\ee
For $Q= Bq$ we have
\be 
U[Q] = R^{\perp} \left ( \fr{1}{1+B^2}Q + W[Q]\right)
\la{UQ}
\ee
with
\be
W[Q] = -\fr{1}{B^2(1+B^2)}\sum_{n=0}^{\infty}\left(\fr{L_Q}{(1+B^2)}\right)^{n+1}Q.
\la{wQ}
\ee
The scalar $w$ and velocity $u$ are unchanged,
\be
W[Q] = w\left [\fr{Q}{B}\right],
\la{WQwq}
\ee
 \be
 U[Q] = u\left [\fr{Q}{B}\right].
 \la{UQuq}
\ee
Let us consider a Banach algebra $\mathcal B$, smaller than $L^{\infty}$ and where the Riesz transforms are bounded. Two canonical examples are
$C^{\alpha}$ and $H^s$, $s> 1$. Then, if $Q\in \mathcal B$ we have that $L_Q$ is bounded in $\mathcal B$ with 
\be
\|L_Q\|_{\mathcal L(\mathcal B, \mathcal B)} \le \gamma \|Q\|_{\mathcal B}.
\la{normLQ}
\ee 
The expressions \eqref{uq}, \eqref{wQ} are explicit analytic expansions in $\mathcal B$, and we have
\be
\|W[Q]\|_{\mathcal B} \le \fr{\gamma}{B^2(1+B^2)^2}\fr{\|Q\|^2_{\mathcal B}}{\left(1 -\fr{\gamma}{1 + B^2}\|Q\|_{\mathcal B}\right)}. 
\la{normwQ}
\ee
If 
\be
\|Q\|_{\mathcal B} \le \fr{1+B^2}{2\gamma}
\la{Qcond}
\ee
then 
\be
\|W(Q)\|_{\mathcal B} \le \fr{2\gamma}{B^2(1+B^2)^2}\|Q\|^2_{\mathcal B}.
\la{normwQB}
\ee
These considerations can be greatly improved by removing the quantitative condition \eqref{Qcond} in $H^s$ spaces. This is done by taking advantage of the unconditional invertibility in $L^2$, provided $Q\in L^{\infty}$, and then using the same cancellation and commutator estimates for higher derivatives.  Because we are interested in the DECM transport equations, we work in a Sobolev space $H^s$ that is smaller than $W^{1, \infty}$ to guarantee Lipschitz velocities, which means that we need $s>2$. For simplicity we take $s$ an integer, $s=3$. 
We have
\beg{prop}\la{goodwQbound}
For any $Q\in H^3$, the equation
\be \label{weq[Q]}
w - \frac{1}{1+B^2} R \cdot (QR^{\perp} w) = - \frac{1}{B^2(1+B^2)^2}R \cdot (QR^{\perp}Q)
\ee
has a unique solution  $w=W[Q]$. There exists an absolute constant $C$ such that $w$ obeys
\be
\|w\|_{L^2}^2 \le \frac{C}{B^4(1+B^2)^4} \|Q\|_{L^4}^4,
\la{w[Q]l2bound}
\ee
and
\be \label{goodw[Q]bound}
\|\na \Delta w\|_{L^2}^2
\le \frac{C(\|Q\|_{L^{\infty}}^5 + \|Q\|_{L^{\infty}}^3 + \|Q\|_{L^{\infty}}^2)}{B^4(1+B^2)^4} \|Q\|_{H^3}^4 + \frac{C\|Q\|_{L^{\infty}}^2}{B^4(1+B^2)^4}\|Q\|_{H^3}^2.
\ee
\end{prop}
We expressed $w$ solving \eqref{weqq} as a function of $Q=Bq$, and thus  rewrote \eqref{weqq}
in terms of $Q$ as \eqref{weq[Q]}.  The function $w$ is the same, $W[Q] = w\left[ \fr{Q}{B}\right]$.
Note that the leading order bound on $\|\na\Delta w\|_{L^2}$ of \eqref{goodw[Q]bound}
relative to the square of the $H^3$ norm is $B^{-2}(1+B^2)^{-2}$, the same as in the abstract bound \eqref{normwQB}, but it  holds without the condition \eqref{Qcond}. The proof of Proposition \ref{goodwQbound} is found in Appendix~\ref{A1}.

We conclude this section by showing some continuity properties of $T_q^{-1}$.
\beg{prop}\label{prop:T}
Let $q \in L^{\infty}$, and let the linear operator $T_{q}: H \rightarrow H$  be defined in \eqref{Tq}
\be 
T_{q}u = (1+B^2) u + B \mathbb P\left (q u^{\perp}\right). 
\ee 
\begin{enumerate}
\item The operator $T_{q}$ is invertible and its inverse $T_{q}^{-1}$ obeys 
\be \label{gal1}
\|T_{q}^{-1} u\|_{L^2} \le \frac{1}{1+B^2} \|u\|_{H}
\ee for all $u \in H$. 
\item If $q_1, q_2 \in L^{\infty}$, then, for every $u \in H$, we have
\be \label{gal}
\|T_{q_1}^{-1} u  - T_{q_2}^{-1} u\|_{H}
\le \frac{B}{(1+B^2)^2} \|q_1 - q_2\|_{L^{\infty}} \|u\|_{H}. 
\ee 
\item Let $\left\{q_n\right\}_{n \in \N}$ be a sequence of functions such that $\left\{q_n\right\}_{n \in \N}$  is uniformly bounded in $L^{\infty}$ and $\left\{q_n\right\}_{n \in \N}$ converges pointwise to $q \in L^{\infty}$  almost everywhere. Then, for every $u \in H$, $\left\{T_{q_n}^{-1} u \right\}_{n \in N}$ converges strongly in $H$ to $T_{q}^{-1} u$. 
\item Let $\left\{q_n\right\}_{n \in \N}$ be a sequence of functions such that $\left\{q_n\right\}_{n \in \N}$  is uniformly bounded in $L^{\infty}$ and $\left\{q_n\right\}_{n \in \N}$ converges strongly in $L^p$ to $q \in L^{\infty}$ for some $p \in (1,\infty)$. Then, for every $u \in H$, $\left\{T_{q_n}^{-1} u\right\}_{n \in N}$ has a subsequence that converges strongly in $H$ to $T_{q}^{-1} u$. 
\end{enumerate}
\end{prop}

\begin{proof} We prove now (1)--(4). 
\begin{enumerate}
    \item The linear bounded operator $J_q(u) = \mathbb P(q u^{\perp})$ is anti-selfadjoint in $H$,
    and therefore $\mathbb I + \fr{B}{1+B^2} J_q$ is invertible with inverse of norm less than 1. Taking the scalar product of the equation  $T_q(u)  = f$ with $u$ and using the anti-symmetry of $J_q$ we have 
    \be 
(1+B^2) \|u\|_{H}^2 = (u,f)_H \le \|u\|_H \|f\|_H
\ee 
and thus $(1+B^2)\|u\|_H \le \|f\|_H$. Because $u = T_q^{-1}f$ we have 
\be
\|T_q^{-1}f\|_H \le \fr{1}{1+B^2} \|f\|_H.
\ee

    \item Let $q_1, q_2 \in L^{\infty}$. We have
    \be 
\|T_{q_1}^{-1}u - T_{q_2}^{-1}u\|_{H}
= \|T_{q_1}^{-1} (T_{q_2} - T_{q_1})T_{q_2}^{-1}u\|_{H}
\le \frac{1}{1+B^2} \|(T_{q_2} - T_{q_1}) T_{q_2}^{-1}u\|_{H}.
    \ee But, for $f=T_{q_2}^{-1}u$,
    \be 
\|T_{q_2}f - T_{q_1}f\|_H = B \|\PP((q_2 - q_1)f^{\perp})\|_H\le B\|q_2-q_1\|_{L^{\infty}}\|f\|_{H}.
\ee  
It follows that 
\be 
\|T_{q_1}^{-1}u - T_{q_2}^{-1}u\|_{H}
\le \frac{B}{1+B^2} \|q_2 - q_1\|_{L^{\infty}}\|T_{q_2}^{-1} u\|_{L^2}
\le \frac{B}{(1+B^2)^2} \|q_1 - q_2\|_{L^{\infty}} \|u\|_{H}
\ee for any $u \in H$.
\item Let $u \in H$. Then we have 
\be 
\beg{aligned}
\|T_{q_n}^{-1} u - T_{q}^{-1}u\|_{H} 
= \|T_{q_n}^{-1}  (T_{q} - T_{q_n})T_{q}^{-1} u\|_{H} 
= \|T_{q_n}^{-1} \left(B \PP ((q - q_n)(T_{q}^{-1}u)^{\perp}) \right)\|_{H}. 
\end{aligned}
\ee Using the uniform boundedness of the operators $T_{q_n}^{-1}$ and the Leray projector $\PP$ on $L^2$, we bound the latter as follows,
\be 
\beg{aligned}
\|T_{q_n}^{-1} u - T_{q}^{-1}u\|_{H} 
\le \frac{B}{1+B^2}  \|(q - q_n) (T_{q}^{-1}u)^{\perp}\|_{L^2}.
\end{aligned}
\ee Since $|(q - q_n) (T_{q}^{-1} u)^{\perp}| \le \left(\|q\|_{L^{\infty}} + \sup\limits_{n \in \N} \|q_n\|_{L^{\infty}} \right)|(T_{q}^{-1}u)^{\perp}|$, $(T_{q}^{-1} u)^{\perp}$ is bounded in $L^2$ by a constant multiple of $\|u\|_{L^2}$, and $q_n$ converges pointwise to $q$ a.e., we deduce that 
\be 
\lim\limits_{n \to \infty} \|(q - q_n) (T_{q}^{-1}u)^{\perp}\|_{L^2} = 0
\ee by the Lebesgue Dominated Convergence Theorem, and consequently 
\be  
\lim\limits_{n \to \infty} \|T_{q_n}^{-1} u - T_{q}^{-1}u\|_{H}  = 0.
\ee 
\item This follows from (3) and the fact that $\left\{q_n\right\}_{n \in \N}$ has a subsequence that converges to $x$ for  a.e. $x \in \TT^2$. 
\end{enumerate}
\end{proof}

\section{Existence of Global Weak Solutions} \label{S2}

In this section, we prove the existence of global weak solutions of \eqref{decmlab} with \eqref{u[q]eq},
for $L^2$ initial charge density on $\TT^2$.  

For each $\epsilon \in (0,1)$, we let $J_{\epsilon}$ be a standard mollifier, and we consider the $\epsilon$-approximate system 
\begin{align}
    &\pa_t q^{\epsilon} + \frac{1}{1+B^2} u^{\epsilon} \cdot \na q^{\epsilon} + \frac{B}{(1+B^2)} J_{\epsilon} R^\perp q^{\epsilon} \cdot \nabla q^{\epsilon} + \frac{1}{(1+B^2)} \Lambda q^{\epsilon} - \epsilon\Delta q^\epsilon= 0, \label{q-epsilon-1}
    \\
    & u^{\epsilon} =  J_{\epsilon} T_{J_{\epsilon}q^{\epsilon}}^{-1} \left[-\PP \left( J_{\epsilon}q^{\epsilon}RJ_{\epsilon}q^{\epsilon} \right)+ B R^{\perp} J_{\epsilon}q^{\epsilon}\right], \la{q-epsilon-2}
    \\
    &q^{\epsilon}(0) = J_{\epsilon} q_0, \la{q-epsilon-3}
\end{align}
in $\TT^2$ with periodic boundary conditions.

\beg{prop}\label{prop:galerkin}
Let $\epsilon \in (0,1)$. Let $q_0 \in L^2$. Then the $\epsilon$-approximate system \eqref{q-epsilon-1}-\eqref{q-epsilon-3} has a unique global smooth solution. 
\end{prop}
The proof of this proposition is found in Appendix~\ref{A2}.

\beg{prop}\label{prop:epsilon-uni-bdd}
Let $\epsilon \in (0,1)$ and $T>0$. Suppose $q_0 \in L^{\infty}$. Then the family $\left\{q^{\epsilon}\right\}$  of solutions of \eqref{q-epsilon-1}-\eqref{q-epsilon-3} is uniformly bounded in $L^{\infty}(0,T; L^{\infty})$ and $L^2(0,T; H^{\frac{1}{2}})$. Moreover, the family $\left\{u^{\epsilon}\right\}$ is uniformly bounded in $L^{\infty}(0,T; L^2)$. 
\end{prop}

\beg{proof}
The time evolution of the $L^2$ norm of $q^{\epsilon}$ is described by the energy balance
\be 
\frac{1}{2}\frac{d}{dt} \|q^{\epsilon}\|_{L^2}^2 + \frac{1}{(1+B^2)} \|\Lambda^{\frac{1}{2}}q^{\epsilon}\|_{L^2}^2 + \epsilon \|\na q^{\epsilon}\|_{L^2}^2 = 0.
\ee
Here we used the fact that $u^\epsilon$ and $\mathcal J_\epsilon R^\perp q^\epsilon$ are divergence-free. Applying the Gr\"onwall inequality, we obtain 
\be
\|q^{\epsilon}(t)\|_{L^2}^2 + \frac{2}{(1+B^2)} \int_{0}^{t} \|\Lambda^{\frac{1}{2}}q^{\epsilon}(s)\|_{L^2}^2 ds
+ 2\epsilon \int_{0}^{t}\|\na q^{\epsilon}(s)\|_{L^2}^2 ds
= \|q^{\epsilon}(0)\|_{L^2}^2 \le \|q(0)\|_{L^2}^2.
\ee for all $t \ge 0$. Now we address the $L^p$ evolution of $q^{\epsilon}$ for even integers $p > 2$. To this end, we multiply the equation obeyed by $q^{\epsilon}$ by $(q^{\epsilon})^{p-1}$, integrate over $\TT^2$, make use of the divergence-free property obeyed by the vector fields $u^{\epsilon}$ and $R^{\perp}q^{\epsilon}$, apply the C\'ordoba–C\'ordoba inequality, and obtain 
\be
    \frac{1}{p}\frac{d}{dt} \|q^{\epsilon}\|_{L^p}^p \leq 0.
\ee
Integrating in time from $0$ to $t$, we infer that 
\be
    \|q^{\epsilon}(t)\|_{L^p} \leq \|q^{\epsilon}(0)\|_{L^p} \le \|q(0)\|_{L^p}
\ee
for any $t \ge 0$. Letting $p\to \infty$, we obtain the $L^{\infty}$ bound \be
\|q^{\epsilon}(t)\|_{L^\infty}\leq \|q(0)\|_{L^\infty}
\ee for all positive times $t \ge 0$. As a consequence of the above estimates, it follows that $\left\{q^{\epsilon}\right\}$ is uniformly bounded in $L^{\infty}(0,T; L^{\infty})$ and $L^2(0,T; H^{\frac{1}{2}})$. Finally, we turn our attention to the regularity of the $\epsilon$-regularized velocities. Due to the boundedness of the operators $\mathcal J_\epsilon$, $T_{\mathcal J_\epsilon q^\epsilon}^{-1}$, $\mathbb P$, and $R$ on $L^2$, we have
\be \label{est:ue}
\begin{aligned}
    \|u^\epsilon\|_{L^2} =& \left\|J_{\epsilon} T_{J_{\epsilon}q^{\epsilon}}^{-1} \left[-\PP \left( J_{\epsilon}q^{\epsilon}RJ_{\epsilon}q^{\epsilon} \right)+ B R^{\perp} J_{\epsilon}q^{\epsilon}\right]\right\|_{L^2} 
    \\
    \leq & C\left\| J_{\epsilon}q^{\epsilon}RJ_{\epsilon}q^{\epsilon}\right\|_{L^2} + C \|q^{\epsilon}\|_{L^2}
    \\
    \leq &C(\|q^\epsilon\|_{L^\infty} + 1) \|q^{\epsilon}\|_{L^2}
\end{aligned}
\ee for some positive constant $C$ depending on $B$ and some universal constants. 
Therefore, the family $\left\{u^{\epsilon}\right\}$ is uniformly bounded in $L^{\infty}(0,T; L^2)$. 

\end{proof}

\beg{prop}\label{prop:epsilon-uni-time-bdd}
Let $\epsilon \in (0,1)$ and $T>0$. Let $q_0 \in L^{\infty}$. The family $\left\{\pa_t q^{\epsilon} \right\}$ 
of solutions of \eqref{q-epsilon-1}-\eqref{q-epsilon-3} is uniformly bounded in $L^2(0,T; H^{-\frac32})$. 
\end{prop}

\beg{proof}
We take an arbitrary test function $\phi \in H^{\frac32}$ and consider the $L^2$ inner product of $\phi$ with each term in \eqref{q-epsilon-1}. For the linear terms, we have
\be
    (\epsilon \Delta q^\epsilon, \phi)_{L^2} \leq \epsilon \|\Lambda^{\frac12} q^\epsilon\|_{L^2} \|\phi\|_{H^{\frac32}}
\ee and 
\be 
 (\Lambda q^\epsilon, \phi)_{L^2} \leq \|q^\epsilon\|_{L^2} \|\phi\|_{H^1}.
\ee
For the nonlinear terms, we  integrate by parts, apply the H\"older and Sobolev inequalities, and obtain 
\be 
    (u^\epsilon\cdot\nabla q^\epsilon, \phi)_{L^2} = - (u^\epsilon\cdot\nabla \phi, q^\epsilon)_{L^2} \leq \|u^\epsilon\|_{L^2} \|\nabla\phi\|_{L^4} \|q^\epsilon\|_{L^4} \leq C\|u^\epsilon\|_{L^2} \|\Lambda^{\frac12} q^\epsilon\|_{L^2} \|\phi\|_{H^{\frac32}}
\ee 
and  
\be
   (J_\epsilon R^\perp q^\epsilon\cdot\nabla q^\epsilon, \phi)_{L^2} = - ( J_\epsilon R^\perp q^\epsilon\cdot\nabla \phi, q^\epsilon)_{L^2}  \leq C\|q^\epsilon\|_{L^2} \|\Lambda^{\frac12} q^\epsilon\|_{L^2} \|\phi\|_{H^{\frac32}}.
\ee 
Putting the above estimates together and using the uniform bounds derived in Proposition~\ref{prop:epsilon-uni-bdd}, we conclude that 
$\left\{\pa_t q^{\epsilon} \right\}$ is uniformly bounded in $L^2(0,T; H^{{-\frac32}})$.

\end{proof}

\beg{thm} \label{wwweak} Let $T>0$, let $q_0 \in L^{\infty}$. Then there exists a global weak solution $q$ 
of  the DECM equation \eqref{decm},  satisfying 
\begin{equation} 
q \in C_{w^*}([0,T]; L^{\infty}) \cap L^2(0,T; H^{\frac{1}{2}}).
\end{equation}   
\end{thm}

\begin{proof}
As $\left\{q^{\epsilon} \right\}_{\epsilon \in (0,1)}$ is uniformly bounded in $L^\infty(0,T;L^\infty)\cap L^2(0,T; H^{{\frac12}})$ and $\left\{\pa_t q^{\epsilon} \right\}_{\epsilon \in (0,1)}$ is uniformly bounded in $L^2(0,T; H^{{-\frac32}})$, there exists a subsequence, also denoted by $\left\{q^{\epsilon} \right\}_{\epsilon \in (0,1)}$, and a limit $q$, such that
\begin{align}\label{conv:q-weak}
    q^\epsilon \stackrel{\ast}{\rightharpoonup} q \text{ in } L^\infty(0,T;L^\infty), \quad q^\epsilon \rightharpoonup q \text{ in } L^2(0,T; H^{{\frac12}}), \quad \partial_t q^\epsilon \rightharpoonup \partial_t q \text{ in } L^2(0,T; H^{{-\frac32}})
\end{align} by the Banach-Alaoglu theorem. 
In addition, by the Aubin-Lions  and Lions-Magenes lemmas, we have
\begin{align}\label{conv:q-strong}
    q^\epsilon \to q \text{ in } L^2(0,T; L^2) \cap C([0,T); H^{-\frac12}).
\end{align}
Since $q\in L^\infty(0,T;L^\infty)$, it follows that $q\in C_{w^*}([0,T); L^\infty)$. 

We denote by 
\begin{equation}\label{eqn:u-limit}
    u := T^{-1}_{q}\left[ -\mathbb P\left(qRq \right) + B R^\perp q  \right].
\end{equation}
and show that $u^\epsilon \to u$ strongly in $L^2(0,T; L^2)$. Denoting
\[
  v^\epsilon := -\PP \left( J_{\epsilon}q^{\epsilon}RJ_{\epsilon}q^{\epsilon} \right)+ B R^{\perp} J_{\epsilon}q^{\epsilon}, \qquad  v := -\mathbb P\left(qRq \right) + B R^\perp q,
\]
and estimating as in \eqref{est:ue}, we have $v^\epsilon, v\in L^\infty(0,T; L^2)$.
We show that $v^\epsilon \to v$ strongly in $L^2(0,T; L^2)$. The convergence of the linear terms follows directly from the boundedness of the Riesz transform on $L^p$ spaces, the convergence property of the mollifiers $\mathcal J_\epsilon$, and the strong convergence of $q^\epsilon$ obtained in \eqref{conv:q-strong}. As for the nonlinear terms, we have
\begin{equation}
    \begin{aligned}
        \|\mathcal J_{\epsilon} q^\epsilon R \mathcal J_{\epsilon} q^\epsilon - qRq\|_{L^2} &\leq \|\mathcal J_{\epsilon} q^\epsilon - q\|_{L^4} \|R \mathcal J_{\epsilon} q^\epsilon\|_{L^4} + \|q\|_{L^\infty} \|R(\mathcal J_{\epsilon} q^\epsilon - q)\|_{L^2}
        \\
        &\leq C \|\mathcal J_{\epsilon} q^\epsilon - q\|_{L^2}^{\frac14} (\|q^\epsilon\|_{L^6} + \| q\|_{L^6})^{\frac34} \|q^\epsilon\|_{L^4} + C\|q\|_{L^\infty} \|\mathcal J_{\epsilon} q^\epsilon - q\|_{L^2}.
    \end{aligned}
\end{equation}
As $q, q^\epsilon\in L^\infty(0,T;L^\infty)$ and $q^\epsilon \to q$ in $L^2(0,T;L^2)$, we deduce that $v^\epsilon \to v$ in $L^2(0,T; L^2)$. We write
\begin{align*}
    u^\epsilon - u &= \mathcal J_\epsilon T^{-1}_{\mathcal J_\epsilon q^\epsilon} v^\epsilon - T^{-1}_q v
    \\
    &= \Big(\mathcal J_\epsilon T^{-1}_{\mathcal J_\epsilon q^\epsilon} (v^\epsilon -  v)\Big) + \Big(\mathcal J_\epsilon (T^{-1}_{\mathcal J_\epsilon q^\epsilon} v - T^{-1}_q v)\Big) + \Big(\mathcal J_\epsilon T^{-1}_q v - T^{-1}_q v\Big):= I_1+I_2+I_3.
\end{align*}
In view of the uniform boundedness of $T^{-1}_{\mathcal J_\epsilon q^\epsilon}$ in $L^2$ and the convergence $v^\epsilon \to v$ in $L^2(0,T; L^2)$, we conclude that $I_1 \to 0$ in $L^2(0,T; L^2)$. Since $\mathcal J_\epsilon q^\epsilon \to q$ in $L^2(0,T; L^2)$ and $v\in L^\infty (0,T; L^2)$, an application of Proposition~\ref{prop:T} yields the convergence of $I_2$ to $0$ in $L^2(0,T; L^2)$. Due to the convergence properties of the mollifiers $\mathcal J_\epsilon$, we have $I_3\to 0$ in $L^2(0,T; L^2)$. Therefore, we conclude that $u^\epsilon \to u$ strongly in $L^2(0,T; L^2)$. 

Finally, let $\phi$ be a test function in $C^\infty([0,T]\times \mathbb T^2).$ By virtue of \eqref{conv:q-weak},  \eqref{conv:q-strong}, and the regularity of $u, u^\epsilon, q, q^\epsilon$, we have 
\[
\langle \partial_t q^\epsilon, \phi\rangle \to \langle \partial_t q, \phi\rangle, \quad \langle \Lambda q^\epsilon, \phi\rangle \to \langle \Lambda q, \phi\rangle, \quad \epsilon \langle \Delta q^\epsilon, \phi\rangle \to 0,
\]
and
\[
  \langle u^\epsilon\cdot\nabla q^\epsilon,\phi\rangle \to \langle u\cdot\nabla q,\phi\rangle, \quad \langle \mathcal J_\epsilon R^\perp q^\epsilon \cdot\nabla q^\epsilon, \phi \rangle \to \langle R^\perp q\cdot\nabla q,\phi\rangle.
\]
By a density argument, it follows that 
\[
\partial_t q + \frac1{1+B^2} u\cdot\nabla q + \frac{B}{(1+B^2)} R^\perp q \cdot\nabla q + \frac{1}{(1+B^2)} \Lambda q = 0 \quad \text{in}\quad L^2(0,T;H^{-\frac32}),
\]
where $u$ obeys \eqref{eqn:u-limit}. The initial condition $q(0)=q_0$ holds due to the weak continuity in time of $q$.

\end{proof}

\section{Convergence to Inviscid SQG} \label{S3}
As mentioned in the introduction, solutions $q(x,\cdot) $ of \eqref{decm} obey \eqref{decmtm} on the first magnetic time scale, that is $q^B(x, t) = q(x, \fr{1+B^2}{B}t)$
obeys
\be \label{inviscid1}
\pa_t q^B + R^{\perp} q^{B} \cdot \na q^B + \frac{1}{B} u^B \cdot \na q^{B} + \frac{1}{B} \Lambda q^B = 0,
\ee 
 with
\be \la{inviscid2}
u^B = - \frac{1}{1+ B^2} \PP\left(q^BRq^B + B q^B(u^B)^{\perp} - {B}R^{\perp}q^B  \right),
\ee in $\TT^2$, with fixed initial data $q^B(0) = q_0$. 

In this section, we prove that any family $\left\{q^B\right\}_{B \ge 1}$ of weak solutions has a subsequence that converges as $B\to\infty $ to a weak solution of the inviscid SQG equation 
\be \la{inviscid3}
\pa_t q + (R^{\perp}q) \cdot \na q = 0
\ee in $\TT^2$. 

\beg{Thm}\la{invsqglim}Let $T> 0$ and $\left\{q^B\right\}_{B \ge 1}$ be a family of weak solutions of the DECM equation \eqref{inviscid1} on $[0,T]$. Then the family has a subsequence that converges weakly as $B\to\infty$ in $L^2(0, T; L^2)$ to a weak solution of the inviscid SQG equation \eqref{inviscid3}.
\end{Thm}

\beg{rem} In the laboratory time scale, any family $q^B$ of weak solutions of \eqref{decm} has a subsequence such that $q^B(x, \fr{1+B^2}{B}t)$ converges weakly as $B\to \infty $ to a weak solution $q(x,t)$ of the inviscid SQG equation \eqref{inviscid3}  in $L^2(0, T; L^2)$.  
\end{rem}

\begin{proof}
Since $\left\{q^{B}\right\}_{B\ge 1}$ is uniformly bounded in $L^2(0,T;L^2)$, there is an unbounded increasing sequence $\left\{B_n\right\}_{n \in \N}$ and a  subsequence of solutions, denoted by $\left\{q^{B_n}\right\}_{n \in \N}$, converging weakly in $L^2(0,T; L^2)$ to some scalar function $q \in L^2(0,T;L^2)$. Due to the uniform boundedness of $\left\{q^{B_n}\right\}_{n \in \N}$ in $L^{\infty}(0,T; L^{\infty})$, the divergence-free property obeyed by $u^{B_n}$ and $R^{\perp} q^{B_n}$, the boundedness of the Riesz transform on $L^2$, and the fact that $B_n \ge 1$, we have 
\be 
\beg{aligned}
\|\pa_t \Lambda^{-1} q^{B_n}\|_{L^2}
&\le \frac{1}{B_n} \|q^{B_n}\|_{L^2} + \frac{1}{B_n} \|\Lambda^{-1} (u^{B_n} \cdot \na q^{B_n})\|_{L^2} + \|\Lambda^{-1} (R^{\perp}q^{B_n} \cdot \na q^{B_n})\|_{L^2}
\\&\le \|q^{B_n}\|_{L^2} + \|u^{B_n}\|_{L^2}\|q^{B_n}\|_{L^{\infty}} + \|R^{\perp} q^{B_n}\|_{L^2}\|q^{B_n}\|_{L^{\infty}}
\\&\le C\left(\|q_0\|_{L^{\infty}}^3 + \|q_0\|_{L^{\infty}}^2 +  \|q_0\|_{L^{\infty}} \right),
\end{aligned}
\ee 
and consequently, we deduce that $\left\{\pa_t \Lambda^{-1}q^{B_n}\right\}_{n \ge 0}$ is uniformly bounded in $L^2(0,T; L^2)$. As $\left\{\Lambda^{-1} q^{B_n}\right\}_{n \in \N}$ is uniformly bounded in $L^2(0,T; H^1)$, there exists a subsequence, also denoted by $\left\{{\Lambda^{-1}}q^{B_n}\right\}$, that converges strongly in $L^2(0,T; H^{\frac{1}{2}})$ to some scalar function $Q$. But $\left\{q^{B_n}\right\}_{n \in \N}$ converges weakly to $q$ in $L^2(0,T;L^2)$, and so $\left\{\Lambda^{-1} q^{B_n}\right\}_{n \in \N}$ converges weakly to $\Lambda^{-1}q$ in $L^2(0,T; {H^1})$. Therefore, $Q = \Lambda^{-1}q$. Now we show that $q$ obeys
\be 
\int_{\TT^2} (q(t) - q_0)\phi(x) dx = -\int_{0}^{t} \int_{\TT^2} (R^{\perp} q \cdot \na q)\phi(x) dx
\ee for a.e. $t \in [0,T]$ and for all $\phi \in H^3(\TT^2)$. To this end, we rewrite \eqref{inviscid1} as 
\be 
q^{B_n}(t) - q_0 = - \int_{0}^{t} R^{\perp} q^{B_n} \cdot \na q^{B_n} ds - \frac{1}{B_n} \int_{0}^{t} u^{B_n} \cdot \na q^{B_n} dx - \frac{1}{B_n} \int_{0}^{t} \Lambda q^{B_n} ds. 
\ee 

As $\left\{\Lambda^{-1} q^{B_n}\right\}_{n \in \N}$ converges strongly in $L^2(0,T; H^{\frac{1}{2}})$, it has a subsequence, still denoted  $\left\{\Lambda^{-1} q^{B_n}\right\}_{n \in \N}$, which converges strongly in $L^2$ to $\Lambda^{-1}q$ for a.e. $t \in [0,T]$, and hence,
\be 
\int_{\TT^2} (q^{B_n} - q)\phi(x) dx = \int_{\TT^2} \Lambda^{-\frac{1}{2}}(q^{B_n} - q) \Lambda^{\frac{1}{2}} \phi dx \rightarrow 0
\ee for a.e. $t \in [0,T]$ and for all $\phi \in H^{\frac{1}{2}}$. In view of the uniform bound $\|q^{B_n}\|_{L^2} \le \|q_0\|_{L^2}$ that holds for all $t \in [0,T]$, we have 
\be 
\beg{aligned}
\frac{1}{B_n}\left|\int_{\TT^2} \int_{0}^{t} \Lambda q^{B_n} \phi dx ds \right|
\le \frac{1}{B_n}\|q_0\|_{L^2} \|\Lambda \phi\|_{L^2}T \to 0, 
\end{aligned}
\ee and thus,
\be 
\frac{1}{B_n} \int_{0}^{t}\int_{\TT^2} \Lambda q^{B_n} \phi dx ds \rightarrow 0 
\ee for all $t \in [0,T]$ and $\phi \in H^1$. By making use of the uniform velocity bound $\|u^{B_n}\|_{L^2} \le C(\|q_0\|_{L^{\infty}}^2 + \|q_0\|_{L^{\infty}})$, we estimate 
\be 
\beg{aligned}
\left|\int_{0}^{t} \int_{\TT^2} u^{B_n} \cdot \na q^{B_n} \phi dx ds \right|
&\le C\int_{0}^{t} \|u^{B_n}\|_{L^2}\|q^{B_n}\|_{L^{\infty}}\|\na \phi\|_{L^2} dx
\\&\le C\|q_0\|_{L^{\infty}}^2(\|q_0\|_{L^{\infty}} +1) \|\na \phi\|_{L^2}T,
\end{aligned}
\ee and hence, 
\be
\frac{1}{B_n} \int_{0}^{t} \int_{\TT^2} u^{B_n} \cdot \na q^{B_n} \phi dx ds \rightarrow 0
\ee for all $t \in [0,T]$ and $\phi \in H^1$. Finally, we consider the difference 
\be 
\int_{0}^{t}\int_{\TT^2} R^{\perp} q^{B_n} \cdot \na q^{B_n} \phi\, dx ds - \int_{0}^{t} \int_{\TT^2} R^{\perp}q \cdot \na q \phi\, dx ds,
\ee and we decompose it into the sum of two terms $\mathcal{I}_n$ and $\mathcal{O}_n$ where 
\be 
\mathcal{I}_n = \int_{0}^{t} \int_{\TT^2} \left(R^{\perp} q \cdot \na (q^{B_n} - q) + R^{\perp}(q^{B_n}-q) \cdot \na q \right)\phi dx ds
\ee and 
\be 
\mathcal{O}_n = \int_{0}^{t} \int_{\TT^2} R^{\perp} (q^{B_n} - q) \cdot \na (q^{B_n} -q)\phi dx ds.
\ee Since $q \in L^{\infty}(0,T; L^{\infty})$, it follows that both quantities $R^{\perp} q \cdot \na \phi$ and $R^{\perp}\cdot (q \na \phi)$ lie in $L^2(0,T; L^{2} )$ for any $\phi \in H^2$, and thus $\mathcal{I}_n \rightarrow 0$ by the weak convergence of $q^{B_n}$ to $q$ in $L^2(0,T;L^2)$. As for the term $\mathcal{O}_n$, we have 
\be 
\beg{aligned}
\mathcal{O}_n 
= -\int_{0}^{t} \int_{\TT^2} R^{\perp} (q^{B_n} - q) (q^{B_n} - q) \cdot \na \phi dx ds
\end{aligned}
\ee via integration by parts, and thus
\be 
\beg{aligned}
\mathcal{O}_n  &= -\int_{0}^{t} \int_{\TT^2} \na^{\perp} \Lambda^{-1} (q^{B_n} - q) (q^{B_n}-q) \cdot \na \phi dx ds
\\&= \int_{0}^{t} \int_{\TT^2} \Lambda^{-1} (q^{B_n} -q) \na^{\perp} (q^{B_n}-q) \cdot \na \phi dx ds
\\&= \int_{0}^{t} \int_{\TT^2} \Lambda^{-1}(q^{B_n} - q) \Lambda \na^{\perp} \Lambda^{-1} (q^{B_n}-q) \cdot \na \phi dxds
\\&= \int_{0}^{t} \int_{\TT^2} R^{\perp} (q^{B_n}-q) \cdot \Lambda \left(\Lambda^{-1}(q^{B_n} - q) \na \phi \right) dxds 
\\&= \int_{0}^{t} \int_{\TT^2} R^{\perp} (q^{B_n}-q) \cdot \left[\Lambda (\Lambda^{-1} (q^{B_n} - q)\na \phi )  - \na \phi \Lambda (\Lambda^{-1}(q^{B_n} - q) \right] dxds - \mathcal{O}_n.
\end{aligned} 
\ee The latter yields the identity 
\be 
\mathcal{O}_n = \frac{1}{2} \int_{0}^{t} \int_{\TT^2} R^{\perp} (q^{B_n} - q) [\Lambda, \na \phi]\Lambda^{-1}(q^{B_n} - q) dxds. 
\ee An application of the Cauchy-Schwarz inequality in the spatial variable gives rise to  
\be 
|\mathcal{O}_n|
\le C \int_{0}^{t} \left(\|q^{B_n} - q\|_{L^2} \|[\Lambda, \na \phi] \Lambda^{-1}(q^{B_n} - q) \|_{L^2} \right) ds.
\ee Using the periodic commutator estimate 
\be 
\|[\Lambda, f]g\|_{L^2}
\le C\|\na f\|_{L^4}\|g\|_{L^4}
\ee that holds for any $f \in W^{1,4}$ and $g \in L^4$, and standard continuous Sobolev embeddings, we bound 
\be 
\|[\Lambda, \na \phi]\Lambda^{-1}(q^{B_n} -q)\|_{L^2} \le C\|\phi\|_{H^3} \|\Lambda^{-\frac{1}{2}} (q^{B_n} - q)\|_{L^2},
\ee and consequently, we obtain 
\be 
|\mathcal{O}_n|
\le C\|\phi\|_{H^3} \|q_0\|_{L^2}\int_{0}^{t} \|\Lambda^{-\frac{1}{2}} (q^{B_n}-q)\|_{L^2} ds \rightarrow 0 
\ee due to the strong convergence of $\Lambda^{-1}q^{B_n}$ to $\Lambda^{-1}q$ in $L^2(0,T; H^{\frac{1}{2}})$. Therefore, $q$ is a weak solution of the inviscid SQG equation.
\end{proof}

\section{Global Regularity and Convergence to Dissipative Critical SQG} \label{S4}

We consider here the DECM equation for the boosted field $Q$. As mentioned in the introduction, if $q(x,\cdot)$ is a solution of \eqref{decm} in laboratory time scale, then 
\be
Q(x,t) = Bq(x, (1+B^2)t)
\la{QqBt}
\ee
solves \eqref{rdecmtm}, which we write as
\be
\pa_tQ + R^{\perp}\left (\left(1+\fr{1}{1+B^2}\right )Q + w\right)\cdot\na Q + \l Q = 0
\la{Qteq}
\ee
with
\be
w = W[Q],
\la{wWQ}
\ee
which is the unique solution of \eqref{weq[Q]} discussed in Proposition \ref{goodwQbound}. We recall that  Proposition \ref{goodwQbound} shows that $w$ is small relative to $\|Q\|^2$ in $H^3$ when $B$ is large. When $B$ is of order 1, this term is difficult to handle and global existence for large data of solutions of the equation is not known, just as in the case of electroconvection in porous media.  We prove here

\beg{Thm}\la{globQ} There exist constants $c>0$, $C>0$ such that if $Q_0\in H^3(\TT^2)$ obeys
\be
\|Q_0\|_{L^{\infty}} \le c
\la{QBcond}
\ee
and 
\be 
\|Q_0\|_{H^3} \le  C B(1+B^2),
\la{Blarge}
\ee
then the solution of \eqref{Qteq}, \eqref{wWQ} with initial data $Q_0$ is unique, exists for all time and obeys
\be
\|Q(t)\|_{H^3} \le \|Q_0\|_{H^3}e^{-\fr{t}{4}}
\la{QH3bound}
\ee
for all $t\ge 0$.
\end{Thm}
\beg{rem}\la{globlab}
In laboratory time scale, the result says that if $q_0$ obeys $B\|q_0\|_{L^{\infty}} \le c$ and $\|q_0\|_{H^3}\le C(1+B^2)$, then
\be
\|q(\cdot,t)\|_{H^3} \le \|q_0\|_{H^3}e^{-\fr{t}{4(1+B^2)}}.
\la{qtlabdecay}
\ee
\end{rem}

\beg{proof}
We start by observing that, without loss of generality we may assume that
\be
\int_{\TT^2}Q(x,t)dx = 0,
\la{intQz}
\ee
because this average is time independent. Secondly, we note that the $L^{\infty}$ norm of $Q$ is nonincreasing in time. We denote by $M$ 
\be
\|Q_0\|_{L^{\infty}} = M
\la{M}
\ee
and we have
\be
\|Q(\cdot, t)\|_{L^{\infty}} \le M
\la{QM}
\ee
a priori. Next, we use the fact that
\be
\| Q\|_{H^3}^2 \sim  \sum_{|\alpha| =3} \int_{\TT^2}|\pa^{\alpha} Q(x,t)|^2dx
\la{normH3}
\ee
because $Q$ has mean zero. Above $\alpha = (\alpha_1,\alpha_2)$ is a multi-index, with $\alpha_j \in \mathbb N$. We recall the
notation $\alpha! = \alpha_1! \alpha_2!$ and $|\alpha| = \alpha_1 + \alpha_2$. 
We compute the evolution 
\be
\fr{1}{2}\fr{d}{dt} \sum_{|\alpha| =3}\int_{\TT^2}|\pa^{\alpha} Q(x,t)|^2dx + D^2 = R_1 + R_2,
\ee
where 
\be
D^2 = \sum_{|\alpha| =3} \|\pa^{\alpha} Q\|_{H^{\fr{1}{2}}}^2,
\la{D}
\ee
and 
\be
R_1 = (1+ \fr{1}{1+B^2})\sum_{|\alpha | =3}\sum_{\beta + \gamma =\alpha, |\beta| >0}\fr{\alpha!}{\beta! \gamma!}\int_{\TT^2}\left(\pa^\beta (R^{\perp}Q)\cdot \na \pa^\gamma Q\right )\pa^{\alpha}Qdx,
\la{a}
\ee  
and 
\be
R_2 =\sum_{|\alpha | =3}\sum_{\beta + \gamma =\alpha, |\beta| >0}\fr{\alpha!}{\beta! \gamma!}\int_{\TT^2} \left(\pa^\beta (R^{\perp}w)\cdot \na \pa^\gamma Q\right )\pa^{\alpha}Qdx.
\la{b}
\ee
The terms with $\beta =0$ vanish because $R^{\perp}Q$ and $R^{\perp}w$ are divergence-free. Therefore $|\gamma| \le 2$. We estimate the two contributions $R_1$ and $R_2$ differently. The first one, like in critical SQG \cite{constantin2001critical}, is going to be absorbed in $D^2$ if $M$ is small,
\be
|R_1| \le C MD^2.
\la{roneneq}
\ee
The second contribution is small compared to $\|Q\|_{H^3}^4$,
\be
|R_2| \le C \fr{1}{B^2(1+ B^2)^2}\|Q\|_{H^3}^4.
\la{rtwoneq}
\ee
We note that 
\be
D^2 \ge C{\|Q\|^2_{H^{3.5}}}
\la{Dlow}
\ee
because $Q$ has mean zero.

For $R_1$ we use $L^3$ bounds.
\be
\|\pa^\beta Q\|_{L^3} \le CM^a D^{1-a}
\la{betaguy}
\ee
where $a =\fr{19 -6|\beta |}{15}$,
\be
\|\na\pa^\gamma Q\|_{L^3} \le CM^b D^{1-b}
\la{gammaguy}
\ee
where $b=\fr{19-6(|\gamma |+1)}{15}$
and
\be
\|\pa^{\alpha} Q\|_{L^3} \le C M^c D^{1-c}
\la{alphaguy}
\ee 
where
$c={\fr{19-18}{15}}$. These are all the same Gagliardo-Nirenberg  inequality
\be
\|\pa^s Q\|_{L^3} \le C\|Q\|_{L^{\infty}}^{\fr{19-6|s|}{15}}\|Q\|_{H^{3.5}}^{\fr{6|s|-4}{15}}
\ee
for $1\le |s|\le 3$.  
Noting that $a+b+c =1$ because $|\beta| + |\gamma| =3$, we proved \eqref{roneneq}.

For the term $R_2$, if $|\beta| =1$ then we estimate $L^{\infty}-L^2-L^2$. 
\be
\left | \int_{\TT^2}(\pa^\beta R^{\perp}w \cdot\na\pa^{\gamma}Q)\pa^{\alpha}Q dx\right |
\le C\|\na R^{\perp}w\|_{L^{\infty}}\|Q\|_{H^3}^2\le C \|w\|_{H^3}\|Q\|_{H^3}^2.
\la{R21}
\ee
If $|\beta| =2$ then we estimate $L^4-L^4-L^2$.
We use the inequality
\be
\|\pa^\beta R^{\perp}w\|_{L^4} \le C\|w\|_{L^2}^{\fr{1}{6}}\|w\|_{H^3}^{\fr{5}{6}}
\la{rbeta}
\ee
\be
\|\na\pa^\gamma Q\|_{L^4} \le C\|Q\|_{L^2}^{\fr{1}{6}}\|Q\|_{H^3}^{\fr{5}{6}}.
\la{nagamma}
\ee
 These are two instances of the inequality
\be
\|\pa^sQ\|_{L^4} \le C \|Q\|_{L^2}^{\fr{1}{6}}\|Q\|_{H^3}^{\fr{5}{6}}
\la{lady}
\ee
for $|s| =2$, so we have for $|\beta|=2$,
\be
\left | \int_{\TT^2}(\pa^\beta R^{\perp}w \cdot\na\pa^{\gamma}Q)\pa^{\alpha}Q dx\right |
\le C\|w\|_{L^2}^{\fr{1}{6}}\|w\|_{H^3}^{\fr{5}{6}}\|Q\|_{H^3}^2\le C \|w\|_{H^3}\|Q\|_{H^3}^2.
\la{R22}
\ee
Finally, for $|\beta|=3$ we estimate using $L^2-L^{\infty}-L^2$ to obtain
\be
\left | \int_{\TT^2}(\pa^\beta R^{\perp}w \cdot\na\pa^{\gamma}Q)\pa^{\alpha}Q dx\right |
\le C\|w\|_{H^3} \|\na Q\|_{L^{\infty}}\|Q\|_{H^3}\le C \|w\|_{H^3}\|Q\|_{H^3}^2.
\la{R23}
\ee
Putting these together we have
\be
|R_2| \le C \|w\|_{H^3}\|Q\|_{H^3}^2.
\la{rtwob}
\ee
Now we use the bound \eqref{goodw[Q]bound} of Proposition \ref{goodwQbound}. Assuming without loss of generality that $M\le 1$ we obtain from \eqref{goodw[Q]bound}
\be
\|w\|_{H^3} \le C \fr{M}{B^2(1+B^2)^2}\|Q\|_{H^3}(1 +\|Q\|_{H^3})
\la{wsmall}
\ee
and using it in \eqref{rtwob} we obtain \eqref{rtwoneq}.

Considering now
\be
y(t) = \sum_{|\alpha| =3}\int_{\TT^2}|\pa^{\alpha} Q(x,t)|^2dx
\ee
and taking $CM\le \fr{1}{2}$ to absorb the term $R_1$ in $\fr{1}{2}D^2$ in view of \eqref{roneneq},
we have
\be
\fr{dy}{dt} + y \le C \fr{1}{B^2(1+B^2)^2} y^2
\la{yineq}
\ee
with some constant $C>0$, where we used a Poincar\'{e} inequality $D^2\ge y$. Thus, if
\be
y(0) \le \fr{1}{3C}B^2(1+B^2)^2
\la{id}
\ee
then, as long as
\be
y(t) \le \fr{1}{2C}B^2(1+B^2)^2,
\la{condy}
\ee
we have that $y$ is decreasing in time and obeys
\be
\fr{dy}{dt} + \fr{1}{2} y \le 0,
\la{odey}
\ee
which implies 
\be
y(t) \le y(0) e^{-\fr{t}{2}}.
\la{ybound}
\ee
But, because $y$ is decreasing, $y(t) \le y(0)$ and because of \eqref{id}, the condition \eqref{condy}
is never violated, so \eqref{ybound} holds for all time. The condition \eqref{Blarge} is the square root of \eqref{id} (with another name for the constant). The bound \eqref{QH3bound} is obtained by taking square roots of both sides of \eqref{ybound}.

\end{proof}

\beg{Thm}\la{critsqgthm} Let $T>0$ and let  $\overline Q$ be a solution on $[0,T]$ of the dissipative critical SQG equation
\be
\pa_t\overline Q + (R^{\perp}\overline{Q})\cdot\na\overline Q + \l \overline Q = 0
\la{critsqg}
\ee
with initial data $\overline Q_0\in H^3$. We consider a strong solution $Q\in L^{\infty}(0,T; H^{\fr{1}{2}}\cap L^{\infty})\cap L^2(0,T; H^1)$ of \eqref{Qteq}, \eqref{wWQ}
with initial data $\overline Q_0$. Then, setting $\overline C = T\exp{\{\int_0^T\|\na \overline Q\|_{L^{\infty}}dt\}}$ we have
\be
\sup\limits_{t\in[0,T]}\|Q(t)-\overline{Q}(t)\|_{L^2} \le \overline C\left (\fr{1}{1+B^2}\|\overline Q_0\|_{L^2}
+ C\fr{1}{B^2(1+B^2)^2}\|\overline Q_0\|_{L^4}^2\right).
\la{l2diff}
\ee
\end{Thm}
\beg{rem} \la{critlab} In laboratory time scale, this result says that strong solutions $q$ of DECM equations \eqref{decm} with initial data $q_0 = \fr{\overline{Q}_0}{B} \in H^3$  obey 
\be
\sup\limits_{t\in[0,T]}\|Bq(\cdot, (1+ B^2)t)-\overline{Q}(t)\|_{L^2} \le\overline C\left (\fr{1}{1+B^2}\|\overline{Q}_0\|_{L^2} + C\fr{1}{B^2(1+B^2)^2}\|\overline{Q}_0\|_{L^4}^2\right) = O\left(\fr{1}{{1+B^2}}\right)
\la{l2q}
\ee
where $\overline Q$ is the global smooth solution of critical SQG \eqref{critsqg} with initial data $\overline{Q}_0 $.
\end{rem}

\beg{proof}
Let $\delta = Q -\overline Q$. 
Writing 
\be
V= \left(1+\fr{1}{1+B^2}\right )Q + w,
\la{VQ}
\ee
we have
\be
\pa_t\delta + \l\delta +(R^{\perp}V)\cdot\na\delta + R^{\perp}( V-\overline Q)\cdot\na\overline Q =0.
\la{deltaeq}
\ee
Because $Q$ is a strong solution we have that $\pa_t Q\in L^{\fr{4}{3}}$, $Q\in L^{\infty}(0,T; L^4)$, and because $\overline Q$ is smooth, we can compute the evolution of $\delta$ in $L^2$, and after cancellations we have 
\be
\fr{1}{2}\fr{d}{dt}\|\delta\|_{L^2}^2 + \|\delta\|_{H^{\fr{1}{2}}}^2  \le \|\na \overline Q\|_{L^{\infty}}\|\delta\|_{L^2}\left( (1+ \fr{1}{1+B^2})\|\delta\|_{L^2} + \fr{1}{1+B^2}\|\overline Q\|_{L^2} + \|w\|_{L^2}\right).
\la{l2ev}
\ee
Using \eqref{w[Q]l2bound} we have
\be
\fr{d}{dt}\|\delta\|_{L^2} \le \|\na \overline Q\|_{L^{\infty}}\left ((1+ \fr{1}{1+B^2})\|\delta\|_{L^2} +  \fr{1}{1+B^2}\|\overline Q_0\|_{L^2}
+ C\fr{1}{B^2(1+B^2)^2}\|\overline Q_0\|_{L^4}^2\right),
\ee
and thus
\be
\sup\limits_{t\in[0,T]}\|\delta(t)\|_{L^2} \le \left (\fr{1}{1+B^2}\|\overline Q_0\|_{L^2}
+ C\fr{1}{B^2(1+B^2)^2}\|\overline Q_0\|_{L^4}^2\right)T\exp{\{\int_0^T\|\na \overline Q\|_{L^{\infty}}dt\}}.
\la{deltal2b}
\ee
\end{proof}

\appendix

\section{Proof of Proposition \ref{goodwQbound}}\label{A1}

{\bf{Step 1. $L^2$ bounds for $w$.}} We multiply the equation \eqref{weq[Q]} by $w$ and integrate over $\TT^2$. Since $R^{\perp} w \cdot R w = 0$, the nonlinear term in $w$ vanishes, and we obtain the identity
\be 
\|w\|_{L^2}^2
= \frac{1}{B^2(1+B^2)^2} \int_{\TT^2} Q R^{\perp} Q \cdot Rw dx. 
\ee Applications of the Cauchy-Schwarz and Young inequalities yield the bound
\be 
\|w\|_{L^2}^2 \le \frac{1}{2} \|w\|_{L^2}^2 + \frac{C}{B^4(1+B^2)^4} \|QR^{\perp}Q\|_{L^2}^2,
\ee from which we deduce that
\be \label{omegal2}
\|w\|_{L^2}^2 \le \frac{C}{B^4(1+B^2)^4} \|Q\|_{L^4}^4,
\ee after making use of the boundedness of the Riesz transform on $L^4$.

{\bf{Step 2. $H^1$ bounds for $w$.}} We take the $L^2$ inner product of the equation \eqref{weq[Q]} obeyed by $w$ with $-\Delta w$. Since $\na$ and $R^{\perp}$ commutes and $R^{\perp} \na w \cdot R \na w =0$, we have  
\be 
\begin{aligned}
\int_{\TT^2} \na R \cdot (Q R^{\perp} w) \cdot \nabla w dx 
= - \int_{\TT^2} \nabla Q R^{\perp} w \cdot R \nabla w dx.
\end{aligned}
\ee Using the interpolation inequality
\be 
\|\na Q\|_{L^{\infty}}
\le C\|Q\|_{L^{\infty}}^{\frac{1}{2}} \|Q\|_{H^3}^{\frac{1}{2}},
\ee we estimate 
\be 
\begin{aligned}
\|\na w\|_{L^2}^2 
&= - \frac{1}{1+B^2} \int_{\TT^2} \na Q R^{\perp} w \cdot R \nabla w dx + \frac{1}{B^2(1+B^2)^2} \int_{\TT^2} \na R \cdot (QR^{\perp}Q) \cdot \na w dx
\\&\le \frac{1}{2} \|\na w\|_{L^2}^2 + \frac{C}{(1+B^2)^2} \|\na Q R^{\perp} w\|_{L^2}^2 
+ \frac{C}{B^4(1+B^2)^4} \|\na R \cdot (QR^{\perp} Q)\|_{L^2}^2
\\&\le \frac{1}{2}\|\na w\|_{L^2}^2 + \frac{C}{(1+B^2)^2} \|Q\|_{L^{\infty}}\|Q\|_{H^3} \|w\|_{L^2}^2 
+ \frac{C}{B^4(1+B^2)^4} \|Q\|_{L^2}^2 \|Q\|_{L^{\infty}} \|Q\|_{H^3}.
\end{aligned}
\ee Due to the boundedness of $w$ in $L^2$ obtained in \eqref{omegal2}, we infer that 
\be \label{omegagrad}
\|\na w\|_{L^2}^2
\le \frac{C(\|Q\|_{L^{\infty}}^5 + \|Q\|_{L^{\infty}}^3)}{B^4 (1+B^2)^4} \|Q\|_{H^3}. 
\ee

{\bf{Step 3. $H^2$ bounds for $w$.}} We apply $\Delta$ to the $w$-equation \eqref{weq[Q]} and take the $L^2$ inner product with $\Delta w$. Using the cancellation law $R^{\perp} \Delta w \cdot R \Delta w = 0$, we have 
\be 
\begin{aligned}
\int_{\TT^2} \Delta R \cdot (Q R^{\perp} w) \cdot \Delta w dx
= - \int_{\TT^2} \left(\Delta Q R^{\perp} w + \na Q \cdot \na R^{\perp} w \right) \cdot R \Delta w dx.
\end{aligned}
\ee Consequently, it holds that 
\be 
\begin{aligned}
\|\Delta w\|_{L^2}^2
&\le \frac{1}{2} \|\Delta w\|_{L^2}^2
+ \frac{C}{(1+B^2)^2} \left(\|\Delta Q\|_{L^4}^2 \|R^{\perp} w\|_{L^4}^2 + \|\na Q\|_{L^{\infty}}^2 \|\na R^{\perp} w\|_{L^2}^2 \right) + \frac{C}{B^4(1+B^2)^4} \|QR^{\perp}Q\|_{H^2}^2
\\&\le \frac{1}{2} \|\Delta w\|_{L^2}^2
+ \frac{C}{(1+B^2)^2} \left(\|Q\|_{L^{\infty}}^{\frac{1}{2}} \|Q\|_{H^3}^{\frac{3}{2}} \|w\|_{L^2}\|\nabla w\|_{L^2} + \|Q\|_{L^{\infty}}\|Q\|_{H^3} \|\na w\|_{L^2}^2 \right)
\\&\quad\quad\quad\quad+ \frac{C}{B^4(1+B^2)^4} \|Q\|_{L^{\infty}}^2 \|Q\|_{H^3}^2
\\&\le \frac{1}{2}\|\Delta w\|_{L^2}^2 +\frac{C(\|Q\|_{L^{\infty}}^5 + \|Q\|_{L^{\infty}}^3 + \|Q\|_{L^{\infty}}^2)}{B^4(1+B^2)^4} \|Q\|_{H^3}^2.
\end{aligned}
\ee where the last two inequalities follow from standard continuous Sobolev embeddings, the boundedness of the Riesz transform on Sobolev and $L^p$ spaces, the fact that $H^2$ is a Banach Algebra, Gagliardo-Nirenberg interpolation inequalities, and application of the bound \eqref{omegagrad}. Thus, we deduce that 
\be \label{omegadelta}
\|\Delta w\|_{L^2}^2
\le \frac{C(\|Q\|_{L^{\infty}}^5 + \|Q\|_{L^{\infty}}^3 + \|Q\|_{L^{\infty}}^2)}{B^4(1+B^2)^4} \|Q\|_{H^3}^2.
\ee 

{\bf{Step 4. $H^3$ bounds for $w$.}} The cancellation law $R^{\perp} \na \Delta w \cdot R \na \Delta w = 0$ gives rise to 
\be 
\begin{aligned}
\|\na \Delta w\|_{L^2}^2
&= -\frac{1}{1+B^2}\int_{\TT^2} \left[\na \Delta (Q R^{\perp} w) - Q R^{\perp} \na \Delta w \right]  \cdot R \na \Delta w dx 
\\&\quad\quad+ \frac{1}{B^2(1+B^2)^2} \int_{\TT^2} \na \Delta (QR^{\perp}Q) \cdot R \na \Delta w dx. 
\end{aligned}
\ee By expanding and simplifying the commutator $[\na \Delta, QR^{\perp}]w$, and estimating using H\"older's inequality and continuous Sobolev embeddings, we obtain 
\be 
\begin{aligned}
\|\na \Delta (Q R^{\perp} w) - Q R^{\perp} \na \Delta w \|_{L^2}^2
\le  C\|Q\|_{H^3}^2 \|\Delta w\|_{L^2}^2.
\end{aligned}
\ee Hence, we have 
\be 
\begin{aligned}
\|\na \Delta w\|_{L^2}^2
\le \frac{1}{2} \|\na \Delta w\|_{L^2}^2 + \frac{C}{(1+B^2)^2} \|Q\|_{H^3}^2 \|\Delta w\|_{L^2}^2
+ \frac{1}{B^4(1+B^2)^4} \|Q\|_{L^{\infty}}^2\|Q\|_{H^3}^2,
\end{aligned}
\ee and so 
\be \label{omegah3}
\|\na \Delta w\|_{L^2}^2
\le \frac{C(\|Q\|_{L^{\infty}}^5 + \|Q\|_{L^{\infty}}^3 + \|Q\|_{L^{\infty}}^2)}{B^4(1+B^2)^4} \|Q\|_{H^3}^4 + \frac{C\|Q\|_{L^{\infty}}^2}{B^4(1+B^2)^4}\|Q\|_{H^3}^2.
\ee due to \eqref{omegadelta}.

\section{Proof of Proposition \ref{prop:galerkin}}\label{A2}

\beg{proof}
For each integer $n \ge 1$, we consider the Galerkin approximants 
\be 
\PP_n \theta = \sum\limits_{j=1}^{n} (\theta, \omega_j)_{L^2} \omega_j
\ee where $\omega_j$ are the eigenfunctions of the negative Laplacian operator with periodic boundary conditions.  For fixed $\epsilon > 0$ and $n \ge 1$, we consider the Galerkin approximate system 
\be \label{galq}
\pa_t q_n^{\epsilon} + \frac{1}{1+B^2} \PP_n (u_n^{\epsilon} \cdot \na q_n^{\epsilon}) + \frac{B}{(1+B^2)} \PP_n(J_{\epsilon} R^\perp q_n^{\epsilon} \cdot \nabla q_n^{\epsilon}) + \frac{1}{(1+B^2)} \Lambda q_n^{\epsilon} - \epsilon\Delta q_n^\epsilon= 0,
\ee 
\be \la{galu}
u_n^{\epsilon} =  \PP_n \left(J_{\epsilon} T_{J_{\epsilon}q_n^{\epsilon}}^{-1}\left[-\PP \left( J_{\epsilon}q_n^{\epsilon}RJ_{\epsilon}q_n^{\epsilon} \right)+ B R^{\perp} J_{\epsilon}q_n^{\epsilon}\right]\right),
\ee with initial data $q_n^{\epsilon} (0) = \PP_n q^{\epsilon}(0)$ and periodic boundary conditions. The latter is a finite-dimensional system of autonomous ODEs and has a unique smooth solution on a maximal time interval. Next, we derive {\it a priori} bounds. 

{\bf{Step 1. $L^2$ bounds.}} We take the $L^2$ inner product of the equation \eqref{galq} obeyed by $q_n^{\epsilon}$ with $q_n^{\epsilon}$. Using the self-adjointness of the Galerkin projectors $\PP_n$, the identity $\PP_n q_n^{\epsilon} = q_n^{\epsilon}$, and the divergence-free condition obeyed by both $u_n^{\epsilon}$ and $R^{\perp} q_n^{\epsilon}$, the nonlinear terms vanish, namely 
\be 
(\PP_{n}(u_n^{\epsilon} \cdot \na q_n^{\epsilon}), q_n^{\epsilon})_{L^2} =  (\PP_n(J_{\epsilon} R^{\perp} q_n^{\epsilon} \cdot \na q_n^{\epsilon}), q_n^{\epsilon})_{L^2}= 0.
\ee This gives rise to the following energy balance
\be 
\frac{1}{2}\frac{d}{dt} \|q_n^{\epsilon}\|_{L^2}^2 + \frac{1}{(1+B^2)} \|\Lambda^{\frac{1}{2}}q_n^{\epsilon}\|_{L^2}^2 + \epsilon \|\na q_n^{\epsilon}\|_{L^2}^2 = 0.
\ee Integrating in time from $0$
to $t$ and using the uniform boundedness of $\PP_n$ on $L^2$, we infer that 
\be \label{gal4}
\|q_n^{\epsilon}(t)\|_{L^2}^2 + \frac{2}{1+B^2} \int_{0}^{t} \|\Lambda^{\frac{1}{2}}q_n^{\epsilon}(s)\|_{L^2}^2 ds
+ 2\epsilon \int_{0}^{t}\|\na q_n^{\epsilon}(s)\|_{L^2}^2 ds
= \|q_n^{\epsilon}(0)\|_{L^2}^2 \le \|q^{\epsilon}(0)\|_{L^2}^2.
\ee for all $t \ge 0$.

{\bf{Step 2. $H^m$ bounds for $m \ge 1$.}} We multiply the equation \eqref{galq} obeyed by $q_n^{\epsilon}$ by $\Lambda^{2m}q_n^{\epsilon}$ and we integrate over $\TT^2$. We obtain the evolution equation 
\be 
\beg{aligned}
&\frac{1}{2} \frac{d}{dt}\|\Lambda^m q_n^{\epsilon}\|_{L^2}^2
+ \frac{1}{(1+B^2)}\|\Lambda^{m+\frac{1}{2}}q_n^{\epsilon}\|_{L^2}^2
+ \epsilon \|\Lambda^{m+1} q_n^{\epsilon}\|_{L^2}^2
\\&= -\frac{1}{1+B^2}\int_{\TT^2} \Lambda^{m-1} (u_n^{\epsilon} \cdot \na q_n^{\epsilon}) \Lambda^{m+1} q_n^{\epsilon} dx 
- \frac{B}{(1+B^2)} \int_{\TT^2} \Lambda^{m-1} (J_{\epsilon} R^{\perp} q_n^{\epsilon} \cdot \na q_n^{\epsilon})\Lambda^{m+1} q_n^{\epsilon} dx.
\end{aligned}
\ee Applications of the Cauchy-Schwarz and Young inequalities give rise to the differential inequality
\be \label{gal3}
\beg{aligned}
\frac{1}{2} \frac{d}{dt}\|\Lambda^m q_n^{\epsilon}\|_{L^2}^2
&+ \frac{1}{(1+B^2)}\|\Lambda^{m+\frac{1}{2}}q_n^{\epsilon}\|_{L^2}^2
+ \frac{\epsilon}{2} \|\Lambda^{m+1} q_n^{\epsilon}\|_{L^2}^2
\\&\le \frac{1}{(1 +B^2)^2} \|\Lambda^{m-1} (u_n^{\epsilon} \cdot \na q_n^{\epsilon})\|_{L^2}^2
+ \frac{B^2}{(1+B^2)^2} \|\Lambda^{m-1} (J_{\epsilon}R^{\perp}q_n^{\epsilon} \cdot \na q_n^{\epsilon})\|_{L^2}^2.
\end{aligned}
\ee By making use of periodic fractional product estimates and standard continuous Sobolev embeddings, we have 
\be \label{nonlin1}
\beg{aligned}
\|\Lambda^{m-1} (u_n^{\epsilon} \cdot \na q_n^{\epsilon})\|_{L^2}^2
&\le C\|\Lambda^{m-1} u_n^{\epsilon}\|_{L^{4}}^2 \|\na q_n^{\epsilon}\|_{L^4}^2
+ C\|u_n^{\epsilon}\|_{L^{\infty}}^2\|\Lambda^{m-1} \na q_n^{\epsilon}\|_{L^2}^2
\\&\le C\|u_n^{\epsilon}\|_{H^{m+1}}^2\|\Lambda^m q_n^{\epsilon}\|_{L^2}^2
\end{aligned}
\ee provided that $m \ge 2$.
We point out that the latter estimate trivially holds when $m=1$ by the H\"older and Sobolev inequalities.
By making use of the explicit relation \eqref{galu}, the boundedness of $\PP_n$ on Sobolev spaces, and the boundedness of the mollifier $J_{\epsilon}$ from $L^2$ to $H^{m+1}$, and the boundedness of $T_{J_{\epsilon}q_n^{\epsilon}}^{-1}$ on $L^2$,  we estimate 
\be 
\beg{aligned}
\|u_n^{\epsilon}\|_{H^{m+1}}^2
&\le \epsilon^{-2m-2} \left\|T_{J_{\epsilon}q_n^{\epsilon}}^{-1}\left[-\PP \left( J_{\epsilon}q_n^{\epsilon}RJ_{\epsilon}q_n^{\epsilon} \right)+ B R^{\perp} J_{\epsilon}q_n^{\epsilon}\right] \right\|_{L^2}^2
\\&\le \epsilon^{-2m-2}(1+B^2)^{-2} \left\|-\PP \left( J_{\epsilon}q_n^{\epsilon}RJ_{\epsilon}q_n^{\epsilon} \right)+ B R^{\perp} J_{\epsilon}q_n^{\epsilon}\right\|_{L^2}^2
\\&\le \epsilon^{-2m-2}(1+B^2)^{-2} \left( \|J_{\epsilon}q_n^{\epsilon}\|_{L^4}^2\|RJ_{\epsilon}q_n^{\epsilon}\|_{L^4}^2 + B^2 \|R^{\perp} J_{\epsilon}q_n^{\epsilon}\|_{L^2}^2 \right).
\end{aligned}
\ee As $J_{\epsilon}$ and $R$ are bounded operators on $L^4$ and $L^2$, and due to the Ladyzhenskaya interpolation inequality, the latter yields
\be 
\beg{aligned}
\|u_n^{\epsilon}\|_{H^{m+1}}^2
&\le C\epsilon^{-2m-2}(1+B^2)^{-2} \left( \|q_n^{\epsilon}\|_{L^4}^4 + B^2 \|q_n^{\epsilon}\|_{L^2}^2 \right)
\\&\le C\epsilon^{-2m-2}(1+B^2)^{-2} \left( \|q_n^{\epsilon}\|_{L^2}^2\|\na q_n^{\epsilon}\|_{L^2}^2 + B^2 \|q_n^{\epsilon}\|_{L^2}^2 \right).
\end{aligned}
\ee Consequently, we obtain the bound
\be 
\|\Lambda^{m-1}(u_n^{\epsilon} \cdot \na q_n^{\epsilon})\|_{L^2}^2
\le C\epsilon^{-2m-2}(1+B^2)^{-2} \left(\|\na q_n^{\epsilon}\|_{L^2}^2  +  B^2\right) \|q_n^{\epsilon}\|_{L^2}^2 
 \|\Lambda^m q_n^{\epsilon}\|_{L^2}^2.
\ee Similarly, we have 
\be 
\beg{aligned}
\|\Lambda^{m-1}(J_{\epsilon}R^{\perp} q_n^{\epsilon} \cdot \na q_n^{\epsilon})\|_{L^2}^2
&\le C\|J_{\epsilon}R^{\perp} q_n^{\epsilon}\|_{H^{m+1}}^2\|\Lambda^m q_n^{\epsilon}\|_{L^2}^2
\\&\le C\epsilon^{-2m-2}\|q_n^{\epsilon}\|_{L^2}^2\|\Lambda^m q_n^{\epsilon}\|_{L^2}^2.
\end{aligned}
\ee Putting all these estimates together, the energy inequality \eqref{gal3} boils down to 
\be
\beg{aligned}
 \frac{d}{dt}\|\Lambda^m q_n^{\epsilon}\|_{L^2}^2
&+ \frac{2}{(1+B^2)}\|\Lambda^{m+\frac{1}{2}}q_n^{\epsilon}\|_{L^2}^2
+ \epsilon \|\Lambda^{m+1} q_n^{\epsilon}\|_{L^2}^2
\\&\le \frac{C}{\epsilon^{2m+2}(1+B^2)^2}\left[\frac{1}{(1 +B^2)^2} \left(\|\na q_n^{\epsilon}\|_{L^2}^2  +  B^2\right) 
+ B^2\right]\|q_n^{\epsilon}\|_{L^2}^2\|\Lambda^m q_n^{\epsilon}\|_{L^2}^2.
\end{aligned}
\ee Finally, we integrate in time from $0$ to $t$, use the uniform bounds \eqref{gal4} derived in Step 1, and deduce that for any $T>0$,
\be 
q_n^{\epsilon} \in L^{\infty}(0,T; H^m(\TT^2)) \cap L^2(0,T; H^{m+1} (\TT^2)).
\ee 

{\bf{Step 3. Convergence.}} The convergence follows from the uniform-in-$n$ boundedness of solutions in the Lebesgue spaces $L^{\infty}(0,T; H^m(\TT^2))$ and $L^2(0,T; H^{m+1} (\TT^2))$ for all $m \in \NN$ and the Aubin-Lions lemma. We point out that the $u_n^{\epsilon}$ converges to $u^{\epsilon}$ due to the uniform boundedness of mollifiers on Sobolev spaces and the Lipschitz estimates \eqref{gal}. The proof is standard and we omit the details. 

{\bf Step 4. Uniqueness.} Suppose there are two smooth solutions $q_1^\epsilon$ and $q_2^\epsilon$ with same initial data and with $u_1^\epsilon$ and $u_2^\epsilon$ determined by $q_1^\epsilon$ and $q_2^\epsilon$ respectively. Denoting the differences by $q^\epsilon=q_1^\epsilon-q_2^\epsilon$ and $u^\epsilon=u_1^\epsilon-u_2^\epsilon$, we have
\be
\beg{aligned}
    &\partial_t q^\epsilon + \frac{1}{1+B^2} (u_1^{\epsilon} \cdot \na q^{\epsilon} + u^{\epsilon} \cdot \na q_2^{\epsilon}) + \frac{B}{(1+B^2)} (J_{\epsilon} R^\perp q_1^{\epsilon} \cdot \nabla q^{\epsilon}+J_{\epsilon} R^\perp q^{\epsilon} \cdot \nabla q_2^{\epsilon})
    \\
    &\qquad + \frac{1}{(1+B^2)} \Lambda q^{\epsilon} - \epsilon\Delta q^\epsilon= 0.
\end{aligned}
\ee 
Multiplying the latter by $q^{\epsilon}$ and integrating over $\TT^2$ give
\be 
\begin{aligned}    &\frac12\frac{d}{dt} \|q^\epsilon\|_{L^2}^2 + \frac{1}{(1+B^2)} \|\Lambda^{\frac12} q^\epsilon\|_{L^2}^2 + \epsilon \|\nabla q^\epsilon\|_{L^2}^2 + \frac{1}{1+B^2} \int_{\mathbb T^2} u^{\epsilon} \cdot \na q_2^{\epsilon} q^\epsilon dx 
    \\
    &\qquad + \frac{B}{(1+B^2)} \int_{\mathbb T^2} J_{\epsilon} R^\perp q^{\epsilon} \cdot \nabla q_2^{\epsilon} q^\epsilon dx = 0.
\end{aligned}
\ee 
Using the boundedness of the operator $T_{J_{\epsilon}q_2^{\epsilon}}^{-1}$ on $L^2$, the Lipschitz estimate \eqref{gal}, and continuous Sobolev embeddings, we estimate $\|u^\epsilon\|_{L^2}$ as follows,
\be 
\begin{aligned}
    &\|u^\epsilon\|_{L^2} = \|u_1^\epsilon-u_2^\epsilon\|_{L^2}
    \\
    \leq & \left\|J_{\epsilon} (T_{J_{\epsilon}q_1^{\epsilon}}^{-1} - T_{J_{\epsilon}q_2^{\epsilon}}^{-1}) \left[-\PP \left( J_{\epsilon}q_1^{\epsilon}RJ_{\epsilon}q_1^{\epsilon} \right)+ B R^{\perp} J_{\epsilon}q_1^{\epsilon}\right] \right\|_{L^2}
    \\
    &\quad\quad\quad\quad+\left\|J_{\epsilon} T_{J_{\epsilon}q_2^{\epsilon}}^{-1} \left[- \PP \left(J_{\epsilon}q_1^{\epsilon}RJ_{\epsilon}q^{\epsilon} + J_{\epsilon}q^{\epsilon}RJ_{\epsilon}q_2^{\epsilon}\right)+ B R^{\perp} J_{\epsilon}q^{\epsilon}\right] \right\|_{L^2}
    \\
    \leq & C \|J_{\epsilon} q^{\epsilon}\|_{L^{\infty}} \left\|-\PP \left( J_{\epsilon}q_1^{\epsilon}RJ_{\epsilon}q_1^{\epsilon} \right)+ B R^{\perp} J_{\epsilon}q_1^{\epsilon} \right\|_{L^2}
    + C \|J_{\epsilon}q_1^{\epsilon}RJ_{\epsilon}q^{\epsilon} + J_{\epsilon}q^{\epsilon}RJ_{\epsilon}q_2^{\epsilon}\|_{L^2} + C \|q^{\epsilon}\|_{L^2} 
    \\
    \leq & \frac C{\epsilon^2} \|q^{\epsilon}\|_{L^2} \left\|-\PP \left( J_{\epsilon}q_1^{\epsilon}RJ_{\epsilon}q_1^{\epsilon} \right)+ B R^{\perp} J_{\epsilon}q_1^{\epsilon} \right\|_{L^2} + C(\|q_1^\epsilon\|_{H^2} + \|q_2^\epsilon\|_{H^2} + 1)  \|q^{\epsilon}\|_{L^2}
    \\
    \leq & \frac C{\epsilon^2} (\|q_1^\epsilon\|_{L^4}^2 + \|q_1^\epsilon\|_{L^2})\|q^{\epsilon}\|_{L^2} + C(\|q_1^\epsilon\|_{H^2} + \|q_2^\epsilon\|_{H^2} + 1)  \|q^{\epsilon}\|_{L^2}.
\end{aligned}
\ee 
Therefore, we obtain the bound
\be
    \int_{\mathbb T^2} u^{\epsilon} \cdot \na q_2^{\epsilon} q^\epsilon dx  \leq C_\epsilon (1+ \|q_1^\epsilon\|_{H^2}^4 + \|q_2^\epsilon\|_{H^2}^2 ) \|q^{\epsilon}\|_{L^2}^2.
\ee 
By making use of the H\"older and Sobolev inequalities, we have
\be
    \int_{\mathbb T^2} J_{\epsilon} R^\perp q^{\epsilon} \cdot \nabla q_2^{\epsilon} q^\epsilon dx \leq C \|q_2^\epsilon\|_{H^2} \|q^{\epsilon}\|_{L^2}^2.
\ee
Combining the above estimates gives rise to the energy inequality
\be
    \frac{d}{dt} \|q^\epsilon\|_{L^2}^2 \leq C_\epsilon (1+ \|q_1^\epsilon\|_{H^2}^4 + \|q_2^\epsilon\|_{H^2}^2 ) \|q^{\epsilon}\|_{L^2}^2.
\ee 
Consequently, the uniqueness of solutions follows from the Gronwall inequality.

\end{proof}

\section*{Acknowledgement}
The work of PC was partially supported by NSF grant DMS-2106528 and by a Simons Collaboration Grant 601960. The work of M.I. was partially supported by NSF grant DMS 2204614.
Q.L. was partially supported by an AMS-Simons travel grant.

\bibliographystyle{plain}
\bibliography{referencev2}

\end{document}